\documentclass[12pt]{article}

\usepackage{fullpage}
\usepackage{comment}
\usepackage{bm}
\usepackage[dvipsnames]{xcolor}
\usepackage{amsmath,amssymb}
\usepackage{amsthm}
\usepackage{wrapfig}
\usepackage{epsfig}
\usepackage{graphicx}
\usepackage{makeidx}
\usepackage{multicol}
\usepackage{tikz}
\usepackage{algorithm}
\usepackage{algpseudocode}
\usepackage{soul}
\usepackage{ulem}
\usepackage{array,multirow,booktabs,bbm} 

\usepackage{algorithmicx}
\usepackage{algpseudocode}
\makeatletter
\algrenewcommand\ALG@beginalgorithmic{\sffamily\small}
\makeatother
\usepackage{epstopdf}
\usepackage{mathrsfs,mathtools}
\usepackage{xcolor}
\usepackage{enumitem}

\newcommand{\N}{{\cal N}}

\graphicspath{{./Figs/}}
\newcommand{\bs}[1]{\boldsymbol{#1}}
\newcommand{\bmu}{\bs{\mu}}

\newcommand{\bc}{\bs{c}}
\newcommand{\bu}{\bs{u}}
\newcommand{\br}{\bs{r}}

\newcommand{\bx}{\bs{x}}

\newcommand{\calP}{\mathcal{P}}

\newcommand{\calF}{\mathcal{F}}

\newcommand{\calN}{\mathcal{N}}

\DeclareMathOperator*{\argmax}{argmax}
\DeclareMathOperator*{\argmin}{argmin}

\theoremstyle{remark}

\begin{document}
\graphicspath{{Figs/}}

\title{AAROC: Reduced Over-Collocation Method with Adaptive Time Partitioning and Adaptive Enrichment for Parametric Time-Dependent Equations

\author{
{Lijie Ji}\thanks{Department of Mathematics, Shanghai University, Shanghai 200444, P.R. China; Newtouch Center for Mathematics of Shanghai University, Shanghai University, Shanghai 200444, P.R. China. Email: {\tt lijieji@shu.edu.cn}. L. Ji acknowledges the support from NSFC grant No. 12201403.},
\quad Zhichao Peng\footnote{Department of Mathematics, The Hong Kong University of Science and Technology, Clear Water Bay, Kowloon, Hong
Kong, China. Email: {\tt  pengzhic@ust.hk}.},\quad
Yanlai Chen\footnote{Department of Mathematics, University of Massachusetts Dartmouth, 285 Old Westport Road, North Dartmouth, MA 02747, USA. Email: {\tt{yanlai.chen@umassd.edu}}. Y. Chen is partially supported by National Science Foundation grant DMS-2208277 and by Air Force Office of Scientific Research, United States grant FA9550-23-1-0037.}}
}

\date{\empty}

\maketitle

\begin{abstract}
Nonlinear and nonaffine terms in parametric partial differential equations can potentially lead to a computational cost of a reduced order model (ROM) that is comparable to the cost of the original full order model (FOM). To address this, the Reduced Residual Reduced Over-Collocation method (R2-ROC) is developed as a hyper-reduction method within the framework of the reduced basis method in the collocation setting. R2-ROC greedily selects two sets of reduced collocation points based on the (generalized) empirical interpolation method for both solution snapshots and residuals, thereby avoiding the computational inefficiency. The vanilla R2-ROC method can face instability when applied to parametric fluid dynamic problems. To address this, an adaptive enrichment strategy has been proposed to stabilize the ROC method. However, this strategy can involve in an excessive number of reduced collocation points, thereby negatively impacting online efficiency.

To ensure both efficiency and accuracy, we propose an adaptive time partitioning and adaptive enrichment strategy-based ROC method (AAROC). The adaptive time partitioning dynamically captures the low-rank structure, necessitating fewer reduced collocation points being sampled in each time segment. Numerical experiments on the parametric viscous Burgers' equation and lid-driven cavity problems demonstrate the efficiency, enhanced stability, and accuracy of the proposed AAROC method.

\end{abstract}

\section{Introduction}

In multi-query problems, such as sensitivity analysis, uncertainty quantification, design optimization, optimal control and inverse problems, multiple simulations of a parametric partial differential equation (pPDE) are required \cite{sinigaglia2022fast,donello2023oblique}. Fast and low-memory simulation tools for pPDEs are highly desired. One attractive approach is the reduced order model (ROM) \cite{HesthavenRozzaStammBook,Quarteroni2015,BennerGugercinWillcox2015,shah2017reduced}. The reduced basis method (RBM) \cite{Rozza2008,Haasdonk2017Review,Quarteroni2015} is a type of projection-based ROM, constructed following an Offline-Online decomposition framework. RBM features a greedy offline stage when a low-dimensional reduced order space is constructed by exploring the low-rank structures of the solution manifold. This iterative algorithm adaptively samples the parameter / configuration space and generate the most informative snapshot, as identified by an error estimator / indicator. In the online stage, computational savings are realized by solving a reduced order problem that projects the full order model (FOM) onto this low-dimensional space generated offline.

For nonlinear nonaffine problems, hyper-reduction techniques \cite{farhat2020computational,carlberg2011efficient} are required to achieve significant computational savings online by avoiding full order evaluation of nonlinear terms. As classified in \cite{farhat2020computational,carlberg2011efficient}, two main types of hyper-reduction techniques are {\it{approximate-then-project}} method and {\it{project-then-approximate}} method. The {\it{approximate-then-project}} approach {first approximates nonlinear nonaffine terms by utilizing a set of reduced collocation points and then project the approximated full order model to a reduced order space.}  
Examples of this approach include Gappy Proper Orthogonal Decompositon \cite{Everson1995}, Empirical Interpolation Method (EIM) \cite{Barrault2004,grepl2005posteriori,drohmann2012reduced,grepl2007efficient}, discrete EIM (DEIM) \cite{chaturantabut2010nonlinear,PeherstorferButnaruWillcoxBungartz2014}, generalized EIM (GEIM) \cite{maday2013generalized}, the best-points interpolation method (BPIM) \cite{nguyen2008efficient,galbally2010non}, missing point interpolation \cite{ryckelynck2005priori} and Gauss–Newton with approximated tensors (GNAT) \cite{CARLBERG2013623,carlberg2011efficient}.
The {\it{project-then-approximate}} approach directly approximates the projected reduced order operators and quantities using reduced collocation points (or reduced meshes) and corresponding weights \cite{legresley2006application}. Hyper-reduction techniques in this category include the adaptive priori hyper-reduction under the finite element framework \cite{ryckelynck2005priori}, Energy Conserving Sampling and Weighting (ECSW) \cite{farhat2014dimensional,klein2023energy}, the Empirical Cubature Method (ECM) \cite{hernandez2017dimensional}, the linear program based Empirical Quadrature Procedure (LP-EQP) \cite{yano2019lp}, and random selected collocation points for nonlinear dynamical systems \cite{bos2004accelerating}.

In this paper, we focus on further developing the Reduced Residual Reduced over Collocation (R2-ROC) method \cite{ChenSigalJiMaday2021}, which can be seen as a specialized \textit{project-then-approximate} approach with identity weights. The reduced collocation points in R2-ROC are sampled using GEIM on solution snapshots and EIM on the full residuals.
Despite the initial success of ROC, directly applying R2-ROC to parametric fluid flow problems may lead to instability in the resulting ROM, 
especially when the dimension of the reduced order space is small. This instability is due to the convection nature of the fluid problem. Intuitively, reduced collocation points should be sampled at locations where the flow has more complex structures. Due to convection effects, such structures may propagate through a large region of the computational domain as time evolves, necessitating a relatively large number of reduced collocation points to capture these structures.   
However, the number of reduced collocation points sampled by the vanilla R2-ROC is always $2n-1$, with $n$ being the dimension of the reduced order space. To improve robustness, the adaptive enrichment strategy-based ROC (AROC) method \cite{chen2022hyper} adaptively samples additional reduced collocation points by monitoring a robustness indicator. Similar over-sampling ideas have also been utilized to improve accuracy and stability in the {\it{approximate-then-project}} framework \cite{Everson1995,peherstorfer2020stability,klein2023energy}. An unintended consequence of the collocation point enrichment strategies is the reduction of the online efficiency.

A key observation is that, although important fluid structures may propagate throughout the entire computational domain, they frequently concentrate locally in small regions for short time periods. As a result, a smaller number of reduced collocation points may suffice for each short time period. To utilize these local structures over short times, we incorporate adaptive time partitioning \cite{Dihlmann2011ModelRO,San2018} and AROC to ensure the robustness and online efficiency, simultaneously. Specifically, we design a robustness indicator for the offline process whose violation will trigger the adaptive enrichment strategy. We then monitor the level of collocation enrichment. If it reaches too high before the robustness criteria is satisfied, the algorithm activates adaptive time partitioning which divides the entire time domain into more segments and restarts the greedy offline process with these newly partitioned time segments. This process will repeat until the stopping criteria of the offline process are met. We call this adaptive time partitioning and adaptive enrichment strategy-based reduced over-collocation (AAROC) method. Via numerical tests of AAROC on Burgers' equation with a discontinuous initial condition and parameterized by a time-dependent viscosity and a lid-driven cavity problem parameterized by the Reynolds number, we observe that the AAROC method achieves the same accuracy as the AROC method with fewer collocation points and is less sensitive to hyperparameters compared to the AROC method.

The paper is organized as follows. The model problem and the full order model are presented in Section \ref{sec:modelandfom}. In Section \ref{sec:aaroc}, we introduce our AAROC method, including adaptive time partitioning and adaptive enrichment strategy for the collocation points. Section \ref{sec:num} presents numerical results for the viscous Burgers' equation and lid-driven cavity problem. Finally, concluding remarks are provided in Section \ref{sec:conclusion}.

\section{Model problem and full order model}
\label{sec:modelandfom}
We consider the following general time-dependent nonlinear parametric partial differential equation (pPDE):
\begin{equation}
\frac{\partial u}{\partial t} - \mathcal{P}(u, t, \bx; \bmu) = g, \quad \bx \in \Omega, \quad t \in [0,T],
\label{eq:ppde}
\end{equation}
where $u(\bx, t)$ denotes the solution to this parametric problem, and $\cal{P}(\cdot, \cdot, \cdot; \bmu)$ is a nonlinear differential operator acting on $u$ and dependent on a parameter $\bmu \in \Omega_p$. Finally, $g(\bx, t)$ is a given forcing term and we assume that boundary and initial conditions are appropriately defined.

The FOM is a high fidelity numerical solver for the model problem \eqref{eq:ppde} associated with the computational grid denoted by $X^\N$. Here, $\N$ represents the degrees of freedom (DoF) in space.  For simplicity, we consider a backward Euler  discretization with a uniform time step size $\Delta t$, resulting in the uniform time nodes given by
\begin{equation}
{\mathcal T}_f = \{t_i=i\Delta t: i = 0, \cdots, \calN_t\}\; \text{ with }\;\calN_t\Delta t=T.
\end{equation}
We assume that the collocation form of the FOM can be expressed as seeking the discrete solution $u(t_{i},X^\N;\bmu)$ associated with the mesh $X^\N$, parameter $\bmu$ and time $t_i$, so that the residual
\begin{align}
&\mathcal{R}(u(t_{i},X^\N;\bmu);u(t_{i-1},X^\N;\bmu))\notag\\
=&\frac{u(t_{i},X^\N;\mu)-u(t_{i-1},X^\N;\bmu)}{\Delta t}-\calP_{\N}(u(t_{i},X^\N;\mu),t;\bmu)-g^\N
\label{eq:full_residual}
\end{align}
equals zero. Here, the discrete differential operator $\calP_{\N}$ approximates $\calP$, and $g^\N$ approximates the source term $g$. Although the full order discrete residual $\mathcal{R}$ is defined with the backward Euler time discretization, our ROM can be straightforwardly integrated with other time discretizations.

\section{ROC method with adaptive time partitioning and adaptive collocation enrichment strategy (AAROC)
\label{sec:aaroc}}
ROC \cite{ChenJiNarayanXu2020,ChenSigalJiMaday2021} follows an Offline-Online decomposition approach. In the offline stage, a low dimensional reduced order space  and a set of reduced collocation points are iteratively constructed in a greedy manner. In the online stage, reduced order solutions defined in the reduced order space are obtained by minimizing the residual only at sampled reduced collocation points. However, like many other hyper-reduced ROMs, ROC may not be robust for certain parametric fluid problems due to the dynamic nature of local structures in such problems. In \cite{chen2022hyper}, an adpative collocation points enrichment strategy is introduced to enhance the robustness of ROC hyper-reduction for nonlinear parametric fluid flow problems. However, this approach may require sampling a relatively large number of reduced collocation points, which can lead to reduced online efficiency.

To ensure the robustness and online efficiency of hyper-reduction method ROC for fluid problems simultaneously, we propose an adaptive time partitioning and adaptive
enrichment strategy-based reduced over-collocation method (AAROC). Adaptive time partitioning \cite{Dihlmann2011ModelRO,cheung2023local,San2015PrincipalID} is introduced to leverage local information in the time structures of the underlying fluid problems, helping to control the number of reduced collocation points.

Before delving into the details, we present the core idea of the offline and online stages of the proposed algorithm. 
\begin{itemize}
    \item[]\textbf{Offline stage:} Adaptively add reduced basis and reduced collocation points in each time segment. If too many reduced collocation points are sampled by the adaptive enrichment strategy, the offline stage will halt for the entire time domain to be split into more segments, and the offline stage will be restarted.
    \item[] \textbf{Online stage:} Obtain the reduced order solution by minimizing the residual at reduced collocation points sampled in each segment.
\end{itemize}

\begin{table}[htbp]
  \begin{center}
  \resizebox{\textwidth}{!}{
    \renewcommand{\tabcolsep}{0.4cm}
    {\scriptsize
    \renewcommand{\arraystretch}{1.5}
    \renewcommand{\tabcolsep}{12pt}
    \begin{tabular}{@{}lp{0.8\textwidth}@{}}
      \toprule
      $\bmu = (\mu_1, \dots, \mu_p)$ & Parameter in $p$-dimensional parameter domain $\Omega_p \subseteq \mathbb{R}^p$ \\
      $\Xi_{\rm{train}}$ & Parameter training set, a finite subset of $\Omega_p$ \\
      $u(\bmu)$ & Function-valued solution of a parameterized PDE on $\Omega_{\bx} \subset \mathbb{R}^{d}$\\
      $\calP(u(\bmu),t,x; \bmu)$ & A (nonlinear) PDE operator\\
 {$N_x$} & Number of finite difference intervals in $x-$ direction of the physical domain\\
      $\mathcal{N}$ & Degrees of freedom (DoF) of a high-fidelity PDE discretization, the ``truth" solver \\ 
      $X^\calN$ & A size-$\calN$ (full) collocation grid\\
      $N_{\rm max}$ & Number of reduced basis snapshots, $N_{\rm max} \ll \mathcal{N}$\\
      $N_{\rm tpar}$ & Number of time segments\\
      $I_{j}$ & The $j$-th time segment, $j=1, \cdots, N_{\text{tpar}}$\\
        $\bmu^j$ & ``Snapshot" parameter values, $j=1, \ldots, N_{\rm max}$\\
        $V_n$ & $n$-dimensional RB space\\
      $\widehat{u}_n(\bmu,t)$ & Reduced basis solution in the $n$-dimensional RB space $V_n$\\
      $E_n(\bmu)$ & Relative Reduced basis solution error \\
      $\Delta_{{N}} \left(\bmu\right)$ & A residual-based error estimate\\
      $X^{n-1}_{r,j}$ & The reduced collocation grid at $n$-th greedy loop in $j$-th segment, a subset of $X^\calN$ determined based on residuals, $j=1, \cdots, N_{\text{tpar}}$\\
      $X^n_{s,j}$ & An additional reduced collocation grid at $n$-th greedy loop in $j$-th segment, a subset of $X^\calN$ determined based on the solutions, $j=1, \cdots, N_{\text{tpar}}$\\
      $X^n_j$ & {The whole reduced collocation grid at $n$-th greedy loop} in $j$-th segment, defined by $X^{n-1}_{r,j} \cup X^{n}_{s,j}$, $j=1, \cdots, N_{\text{tpar}}$\\
      $T$ & Final time for the time-dependent problems\\
      $\Delta t$ & Time stepsize for the time dependent problems\\
        $t_i$ & Time level $i$, $i=1, \ldots, \calN_t$\\
          $\calN_t$ & Total number of time levels, i.e. $\calN_t = T/\Delta t$\\
      ${\mathcal T}_f$ & {The set of uniformly discretized time nodes $t_i$ of the entire time domain $[0,T]$} \\
        ${\mathcal T}_j$ & {The set of uniformly discretized time nodes $t_i$ within $j$-th time segment $I_j$, $j=1, \cdots, N_{\text{tpar}}$}\\
      ${\mathcal T}_r$ & The set of selected time nodes for the RB space during the offline process\\
        $\rho_\Delta^{n}$ &  The robustness indicator of $n-$th greedy loop\\
      $\gamma$ & Control parameter for the robustness indicator of the AAROC method\\
      $p_{\text{adap}}$ & Percent parameter for the adaptive enrichment strategy\\
      $n_{\text{adap}}^{\rm{incre}}$ & Stepsize of the adaptive enrichment strategy\\
      $n_{\rm{adap}}^{\max}$ & The maximum number of adaptive enriched collocation points \\
      $\epsilon_{\mathrm{tol}}$ & Error estimate stopping tolerance in greedy sweep \\
      \midrule
      Offline process & The pre-computation phase, during which the reduced basis and reduced collocation points in each time segment are determined through a greedy loop.
\\
      Online process & The process of solving the reduced problem to obtain the reduced order solution.\\
    \bottomrule
    \end{tabular}
  }
    }
  \end{center}
\caption{Notation and terminology used throughout this article.}\label{table:notation}
\end{table}
For the convenience of readers, we summarize commonly used notations in Table \ref{table:notation}.

\begin{figure}[htb]
    \centering
    \includegraphics[width=0.95\linewidth]{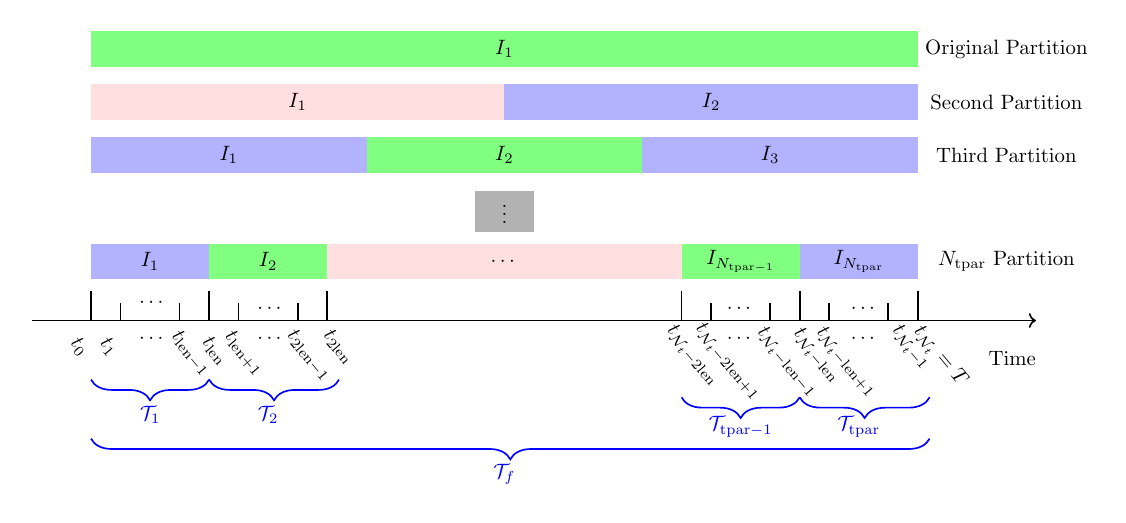}
    \caption{Schematic illustration of the adaptive time partitioning technique of AAROC method. Here, $\text{len}={\mathcal{N}_t}/{N_{\text{tpar}}}$, and ${\mathcal T}_j$ is the set of uniformly discretized time nodes $t_i$ within $j$-th time segment $I_j$, $j=1, \cdots, N_{\text{tpar}}$.}
    \label{fig:schfigure}
\end{figure}

\subsection{Online stage}
In the offline stage, we construct an $n$-dimensional RB space $V_n$, partition the time domain into  $N_{\text{tpar}}$ successive segments $I_j$, $j=1, \cdots, N_{\text{tpar}}$, and {select a set of} reduced collocation points $X_{j}^m \subset X^\N$ {for} each segment $I_j$, {see Figure \ref{fig:schfigure}}. The corresponding collocation operator $P_{*,j}$ for each segment $I_j$ can be defined accordingly. 

In the online stage, given the reduced order solution at time $t_{i-1}$, namely $\widehat{u}_n(\bmu, t_{i-1})$, we seek the reduced order solution at any time $t_i = i\Delta t$ for parameter $\bmu$, 
$\widehat{u}_n(\bmu, t_i) = V_{n} \bc_n (\bmu, t_i)$,
 by minimizing the residual defined in \eqref{eq:full_residual}, constrained at sampled reduced collocation points:
\begin{align}
\bc_n(\bmu, t_i) =
\argmin_{\omega}\left\lVert
P_{*,j}\left(\mathcal{R}(V_n \omega;\widehat{u}_n(\bmu, t_{i-1}))\right)
\right\rVert.
\label{eq:pde:reduced_t}
\end{align}
Here, {the collocation operator $P_{*,j}$ is determined by the time $t_i$,} and the residual is evaluated only at the reduced collocation points, eliminating the need for full order evaluations. 
In practice, evaluating the reduced residual can be achieved by calculating specific rows from the discretized matrices of the nonlinear and linear operators, as well as the discretized forcing term.

\subsection{Offline stage}
In the offline stage, the proposed AAROC method iteratively updates the reduced basis and samples reduced collocation points using greedy algorithms. The key difference between R2-ROC and standard Galerkin RBM is the introduction of collocation points and an error estimator built on the reduced collocation girds. To improve the stability of R2-ROC, the AROC method further introduces an adaptive enrichment strategy for the collocation points. Despite its enhanced stability, the online efficiency of AROC may be reduced when too many collocation points are sampled. To control the number of collocation points, we introduce an adaptive time partitioning strategy.

In Secs. \ref{sec:rb-basis} to \ref{sec:time-partition}, we present essential details of each component of the proposed AAROC method and summarize the algorithm in Alg. \ref{alg:aaroc}. Throughout Secs. \ref{sec:rb-basis} to \ref{sec:time-partition}, we assume that the time domain $[0,T]$ is partitioned into segments $\cup_{j=1}^{N_{tpar}}I_j$.  We denote the time nodes inside $I_j$ as $\mathcal{T}_j$ and the set of reduced collocation points associated with $I_j$ and the $n-$dimensional RB space as $X_j^n$.

\subsubsection{Greedy basis generation and reduced error indicator\label{sec:rb-basis}}

The reduced space is expanded by adding solutions corresponding to the ``underrepresented" $(\bmu,t)$-pair in each greedy iteration. ROC method first samples the ``most underrepresented" parameter $\bmu$ identified by an error indicator, and then samples the ``most underrepresented" time $t$ corresponding to this selected $\bmu$. Under the R2-ROC framework, a reduced error indicator defined on sampled reduced collocation points is used. The detailed procedure is outlined below.

\textbf{Initialization:} The first parameter {$\bmu^1$} is randomly sampled. The first time $t_{\bmu^1}^{1}$ is determined by maximizing the variation across all time instants:
\[
  t_{\bmu^1}^{1} = \argmax_{t_i \in {\mathcal T}_f} \left( \max_{x \in X^\N} \bu(t_i, X^\N; \bmu^1)- \min_{x \in X^\N} \bu(t_i, X^\N; \bmu^1) \right),
\]
where $\{\bu(t_i, X^\N; \bmu^1)\}_{i=0}^{\calN_t}$ are snapshots provided by the high fidelity full order model. {The RB space $V_1$ is initialized as $V_1 =\xi_{1}= \{\bu(t_{\bmu^1}^{1}, X^\N; \bmu^1)\}$.}

\textbf{Reduced error indicator:}
For each $\bmu$ in the training set $\Xi_{\rm train}$ a discretization of the parameter domain $\Omega_p$, we compute the RB approximation $\widehat{\bu}_{n}(\boldsymbol{\mu},t_i)=V_n \bc(\mu,t_i)$ using a procedure similar to the online stage. 
The full residual for the reduced order solution is defined as 
\begin{equation}
\br_n(t_i; \bmu) = V_n \frac{\bc(t_i;\bmu)-\bc(t_{i-1};\bmu)}{\Delta t} -\calP_\N(V_n\bc(t_i;\bmu), t_i; \bmu) -g^\N.\label{eq:fullresidual:rb}
\end{equation}
The error indicator for the entire time history is defined as 
\begin{align}
\Delta_{n}^{RR_t}(\bmu) \coloneqq \sum_{j=1}^{N_{\rm{tpar}}}  \sum_{t_i \in {\mathcal T}_{j}}  \varepsilon^{RR}_j(t_i; \bmu).\label{eq:tpar2:delta_t}
\end{align}
where 
\begin{align}
&\varepsilon^{RR}_j(t_i; \bmu) \coloneqq \lVert \br_{n,j}(t_i; \bmu)\rVert_\infty,\quad \text{if }t_i\in I_j,\quad j=1, \cdots, N_{\text{tpar}},\\
&\br_{n,j}(t_i;\bmu) = P_{\ast,j}\br_n(t_i;\bmu). \label{eq:reduced-residual}
\end{align}
Here, $P_{*,j}$ is the collocation operator associated with the time segment $I_j$, and $\br_{n,j}(t_i;\mu)$ is the reduced residual evaluated only at reduced collocation points. In other words, to compute the error indicator, we evaluate the residual only at these collocation points.

{\bf Greedy in $\bmu$:} The next parameter $\bmu$ is sampled  by maximizing the error indicator:
\begin{align}
  \bmu^{n+1} = \argmax_{\bmu \in \Xi _{\rm train}} \Delta_n^{RR_t}(\bmu),~~
  \Delta_n = \Delta_n^{RR_t}(\bmu^{n+1}).
  \label{eq:errorestimator:rb}
\end{align}

{\bf Greedy in $t$:} 
Once a new parameter $\bmu^{n+1}$ is determined, we compute the RB approximations $\widehat{\bu}_n(\bmu^{n+1}, t) = V_{n} \bc_n (\bmu^{n+1}, t)$ at all time levels $t \in {\mathcal T}_f$. Then, a new time sample is greedily selected by maximizing the reduced residuals at all time levels: 
\begin{align}
t_{\bmu^{n+1}}^{k_{\bmu^{n+1}}} \coloneqq \argmax_{t \in {\mathcal T}_f}\left\{\varepsilon^{RR}( t;\bmu)\right\}, \mbox{ and } {\mathcal T}_r = {\mathcal T}_r \bigcup \{t_{\bmu^{n+1}}^{k_{\bmu^{n+1}}}\}.
 \label{eq:errorestimator:rbt}
\end{align}
Here, ${\mathcal T}_r$ represents the set of sampled time nodes. The corresponding high fidelity solution 
 $u(t_{\bmu^{n+1}}^{k_{\bmu^{n+1}}},X^\calN;\bmu^{n+1})$ is the $(n+1)$-th new basis, denoted by $\xi_{n+1}$.

\subsubsection{ROC collocation points}
After enriching the reduced space in each greedy iteration, we first add two collocation points in each segment $I_j$ following the ``double greedy procedure" of the R2-ROC framework \cite{ChenSigalJiMaday2021}. The key is to perform the ``double greedy procedure" separately within each segment.

\textbf{GEIM collocation points for newly added basis function:}
Define the new basis function resulting from adding the new greedy sample $(\bmu^{n+1},t_{\bmu^{n+1}}^{k_{\bmu^{n+1}}})$ as $\xi_{n+1}$. A collocation point is added using the GEIM method \cite{MadayMulaPateraYano2015}, with the GEIM interpolating functional defined as:
\begin{equation}
\sigma_{\bx, t}^{\bmu}(v) = \frac{v(t,X^\N;\bmu)}{\Delta t}-\mathcal{P}_{\N}(v(t,X^\N;\mu),t;\mu),\quad t\in [0,T].
\label{eq:geimfunctional}
\end{equation}
Using this interpolating functional, GEIM ensures that the discrete differential operators applied to the reduced order solution and the full order solution yield the same value at the interpolation points. 
For detailed procedures on how to add collocation points using GEIM, please refer to Alg. \ref{alg:GEIM}.

\textbf{EIM collocation points for the residual:}
In each time segment $I_j$, we greedily sample an underrepresented $(\bmu,t)$ pair as follows:
\begin{align}
\begin{cases}
 & \bmu^{n+1}_j = \argmax_{\bmu \in \Xi _{\rm train}} \Delta_{n,j}^{RR_t}(\bmu),\\
&t_{\bmu_j^{n+1}}^{k_{\bmu_j^{n+1}}} \coloneqq \argmax_{t \in {\mathcal T}_j}\left\{\varepsilon^{RR}( t;\bmu_{j}^{n+1})\right\}.
\end{cases}
  \label{eq:errorestimator:rb2}
\end{align}
We then apply the EIM to the full residual $\br_{n}(t_{\bmu^{n+1}_j}^{k_{\bmu^{n+1}_j}}; \bmu^{n+1}_j)$ to add additional collocation points (see Alg. \ref{alg:EIM}).

\begin{algorithm}[h]
\begin{algorithmic}[1]
\vspace{0.5ex}
\State \textbf{Input:} interpolation operators $\{\sigma_{1,j}, \cdots, \sigma_{n,j}\}$, collocation set $X^n_j$, sub-sampling operator matrix $P_{*,j}$, and new basis $\xi_{n+1}$
\State Compute a generalized interpolatory residual for vector $\xi_{n+1}$:

Find $\xi_{n+1} = {\xi_{n+1}} - V_{n}\alpha$ with $\{\alpha_i\}$ determined by {requiring that} $\sigma_{i,j}(\xi_{n+1})=0$ for $i= 1, \ldots, n$

Find $\bx^{n+1}_{*,j}=\argmax_{\bx \in X^\N/X^n_j} \left|\sigma_{\bx, t_{\bmu^{n+1}}^{k_{\bmu^{n+1}}}}^{\bmu^{n+1}}(\xi_{n+1})\right|$

Define $\sigma_{n+1,j}(\cdot) := \sigma_{\bx^{n+1}_{*,j}, t_{\bmu^{n+1}}^{k_{\bmu^{n+1}}}}^{\bmu^{n+1}}(\cdot)$. Let $i_{1}$ be the index of $\bx^{n+1}_{*,j}$ in $X^\N$ and update the 

sub-sampling operator matrix $P_{*,j} = [P_{*,j}; \,\, (e_{i_1})^T]$

\State \textbf{Output:} $\{\sigma_{1,j}, \cdots, \sigma_{n+1.j}\}$, $\{\bx^{n+1}_{*,j}\}$, $P_{*,j}$
\end{algorithmic}
\caption{\sf GEIM algorithm in segment $I_j$ for the snapshots}
\label{alg:GEIM}
\end{algorithm}

\begin{algorithm}[h]
\begin{algorithmic}[1]
\vspace{0.5ex}
\State \textbf{Input:} sub-sampling operator matrix $P_{*,j}$, residual vectors $\{r_{1,j}, \cdots, r_{n-1,j}\}$, collocation set $X^n_j$, residual-based collocation set $X_{r,j}^{n-1}$, collocation $\bx_{*,j}^{n+1}$ and parameter-time pair $(\bmu^{n+1}_j, t_{\bmu^{n+1}_j}^{k_{\bmu^{n+1}_j}})$
\State Calculate the full residual $r_{n}(t_{\bmu^{n+1}_j}^{k_{\bmu^{n+1}_j}}, \bmu^{n+1}_j)$ by Eq.~\eqref{eq:fullresidual:rb}

\State Compute an interpolatory residual $r_{n,j}$: 

Find $r_{n,j} = r_{n,j} - \sum_{k=1}^{n-1}\alpha_k r_{k,j}$ with $\{\alpha_k\}$ determined by {requiring that} $r_{n,j}(X^{n-1}_{r,j})=0$

Find $\bx_{**,j}^n=\argmax_{\bx \in X^\N/\left\{X^n_j, \bx_{*,j}^{n+1}\right\}} \lvert r_{n,j}\rvert$

Calculate $r_{n,j}=r_{n,j}/ r_{n,j}(\bx_{**,j}^n)$. Let $i_2$ be the index of $\bx_{**,j}^n$ in $X^\N$ and update the 

sub-sampling operator matrix $P_{*,j} = [P_{*,j}; \,\, (e_{i_2})^T]$

\State \textbf{Output:} $\{r_{1,j}, \cdots, r_{n-1,j}, {r_{n,j}}\}$, $\{\bx^{n}_{**,j}\}$, $P_{*,j}$. 
\end{algorithmic}
\caption{\sf EIM algorithm in segment $I_j$ for the residuals}
\label{alg:EIM}
\end{algorithm}

\subsubsection{Adaptive time partitioning and adaptive collocation enrichment strategy}
\label{sec:time-partition}
As demonstrated in the AROC method \cite{chen2022hyper}, simply adding two reduced collocation points in each greedy iteration may lead to instability in parametric fluid problems. 
Here, we apply the adaptive enrichment strategy from the AROC method independently in each time segment $I_j$. To avoid significant reduction in online efficiency due to an excessive number of collocation points, we further adaptively partition the time domain. 

\textbf{Robustness indicator:}
A {\textit {robustness indicator}} $\rho_\Delta^{n}$ is defined as the ratio of two adjacent error estimators once a sufficient number of reduced basis has been generated
\[
\rho_\Delta^{n}=\begin{cases}
\frac{\Delta_n}{\Delta_{n-1}}, \quad n\geq n_0, \\
0,\quad n<n_0.
\end{cases}
\]

\textbf{Adaptive time partitioning:} We adaptively augment the reduced collocation set in each time segment $I_j$ until the robustness indicator is smaller than a selected tolerance or the maximum number of collocation points allowed in one segment is reached. If the maximum number of collocation points allowed in each segment is reached before the robustness criterion is satisfied, we increase the number of time segments by one through uniform partitioning of the time domain. 

\textbf{Adaptive collocation enrichment strategy:} when the robustness indicator exceeds a tolerance, we revisit the full residual from the previous iteration $r_{n-1,j}$ and choose additional collocation points to enhance the $n$-dimensional ROC solver and regenerate the new reduced basis.

Following \cite{chen2022hyper}, we use an inverse cumulative distribution function (CDF), which is defined as the percentile of $\bx$ within the discretized sorted full collocation set, to adaptively add more collocation points. 
Specifically, we calculate the CDF $\calF_{\lvert\br_{n-1,j}\rvert}$ of $\br_{n-1,j}(\bx)$ across all s by sorting the full order residual $\br_{n-1,j}(\bx)$ in decreasing order of its absolute values. 
We sample an additional $n_{\rm{adap}}$ collocation points uniformly from the top $p_{\rm{adap}}$ (e.g. $20\%$) of the sorted set. These additional points in each segment are defined as: 
\begin{equation}
\bx_{+,j}^{n-1} = \calF_{\lvert\br_{n-1,j}\rvert}^{-1}\left (1-p_{\rm{adap}}: \frac{p_{\rm{adap}}}{n_{\rm{adap}} - 1}: 1\right).
\label{eq:add_colo_formula}
\end{equation}
Here, $a:h:b$ is the Matlab notation for sampling the interval $[a,b]$ with a step size $h$. Consequently,
the second set of residual-based collocation points is updated as $X^{n-1}_{r,j} =X^{n-1}_{r,j} \cup \{\bx_{+,j}^{n-1}\}$, and the entire collocation set is updated as $X^n_j =X^n_j \cup \{\bx_{+,j}^{n-1}\}$. 
The error indicator $\Delta_n$ and the {\textit{robustness indicator}} $\rho_\Delta^{n}$ are then updated. A new $(n+1)$th parameter-time pair will be determined.

\textbf{AAROC algorithm:} Integrating all components, we now have the AAROC algorithm. We summarize the adaptive enrichment strategy in 
Alg. \ref{alg:adaptiveTimePartition} and the complete AAROC algorithm in Alg. \ref{alg:aaroc}.
A key difference between AAROC and AROC is that AAROC avoids oversampling by AROC by the introduction of time partitioning (i.e., the outer while loop in Alg. \ref{alg:aaroc}).  
This mechanism allows AAROC to ensure both accuracy and efficiency by controlling the number of reduced collocation points sampled in each time segment (i.e., the inner adaptive loop Alg. \ref{alg:aaroc}).

\begin{algorithm}[h]
\begin{algorithmic}[1]
\vspace{0.5ex}
\State \textbf{Input:} control parameters $\gamma$, error estimators $\Delta_{n-1}$ and $\Delta_{n}$, full residual $r_{n-1,j}$, residual-based 

collocation set $X_{r,j}^{n-1}$, 
solution-based collocation set $X^{n}_{s,j}$, collocation set $X^n_j$, interpolatory

operator matrix $P_{*,j}$, $n_{\text{adap}}^{\rm{incre}}$, $n_{\rm{adap}}^{\max}$. Initialize $n_{\text{adap}} =n_{\text{add}}$.
\State \mbox{\textbf{While}} $\rho_\Delta^{n} > \gamma \& n_{\rm{adap}} \leq n_{\rm{adap}}^{\max}$
\State Revisit full residuals $r_{n-1,j}$ in all segments and sample $n_{\text{adap}}$ elements $\bx_{+,j}^{n}$ from the top $p_{\text{adap}}$ 

percentile of the domain of $r_{n-1,j}$. Update collocation sets, solve the reduced system, 

Eq. \eqref{eq:pde:reduced_t}, and re-calculate indicator value $\rho_\Delta^{n}$.
\State $n_{\rm{adap}} = n_{\rm{adap}}+ n_{\rm{adap}}^{\rm{incre}}.$
\State \mbox{\textbf{End While}} 
\State Update collocation sets $X^{n-1}_{r,j} = X^{n-1}_{r,j} \cup \bx_{+,j}^{n}$, 
$X^n_j = X^{n}_{s,j} \cup X_{r,j}^{n-1}$, operator $P_{*,j}=[P_{*,j};e_{i_+}^T]$.
\State \mbox{\textbf{If}} $\rho_\Delta^{n} > \gamma$, then $\textrm{I}_{\textrm{aa}}=0$, otherwise $\textrm{I}_{\textrm{aa}} =1$.
\State \textbf{Output:} $X^{n}_{s,j}$, $X^{n-1}_{r,j}$, $X^n_j$, $P_{*,j}$, $\Delta_{n}$, $\textrm{I}_{\textrm{aa}}$
\end{algorithmic}
\caption{\sf Adaptive collocation enrichment algorithm}
\label{alg:adaptiveTimePartition}
\end{algorithm}

\begin{algorithm}[h]
\begin{algorithmic}[1]
\vspace{0.5ex}
\State \textbf{Input:} control parameters $p_{\text{adap}}$, $n_{\text{add}}$, $n_{\text{adap}}^{\rm{incre}}$, $n_{\rm{adap}}^{\max}$, $N_{\max}$, $N_{\text{tpar}}$, $\gamma$, training set $\Xi_{\rm train}$, number

of time partitioning $n_{\text{tpar}}=1$, $n=1$
\State \mbox{\textbf{While}} {$n_{\text{tpar}} < N_{\text{tpar}}$}
\State Select the first parameter-time pair $(\bmu^1,t_{\bmu^{1}}^{k_{\bmu^{1}}})$, and generate first basis $\xi_1$ and first collocation

point $x^\ast$.
\State  Initialize $n=1$, $n_{\rm{adap}}^{\max}=n_{\rm{add}}$, solution-based collocation set $X^{n}_{s,j}$, interpolatory operator  matrix

 $P_{\ast,j}$, $j=1, \cdots, n_{\text{tpar}}$.
\State \mbox{\textbf{While}} {$n < N_{max}$}
\State{Solve the reduced system \eqref{eq:pde:reduced_t} to obtain $\bc_{n}(\bmu, t_k)$, then calculate $\varepsilon^{RR}_j (t_k, \bmu)$ 
and sub-error indicator

 $\Delta_{n,j}^{RR_t}(\bmu)$, $j=1, \cdots, n_{\text{tpar}}$.}
\State Calculate error indicator by Eq.\eqref{eq:tpar2:delta_t}.
\State Perform Algorithm \ref{alg:adaptiveTimePartition}. If  $\textrm{I}_{\textrm{aa}} =0$, break, otherwise, select the $(n+1)$th parameter-time pair by Eqs.

  \eqref{eq:errorestimator:rb} and \eqref{eq:errorestimator:rbt}. 
Obtain the $(n+1)$th basis $\xi_{n+1}$ through the high fidelity solver and generate a 

collocation point $\bx_{\ast,j}^{n+1}$ for each segment through Algorithm \ref{alg:GEIM}.
\State  Determine another parameter-time pair for each segment by Eq. \eqref{eq:errorestimator:rb2} and generate second 

collocation point $\bx_{**,j}^{n}$ 
for each segment through Algorithm \ref{alg:EIM}.
\State  Update collocation sets $X^{n}_{r,j} = X^{n-1}_{r,j} \cup \bx_{\ast \ast,j}^n$, $X^{n+1}_{s,j} = X^{n}_{s,j} \cup \bx_{\ast,j}^{n+1}$, 
$X^{n+1}_j = X^{n+1}_{s,j} \cup X_{r,j}^{n}$, 

$m= \#\{X^{n+1}_{s,j},X^{n}_{r,j}\}$, 
basis space $V_{n+1} =\{V_n, \xi_{n+1}\}$, $n=n+1$, and the operator matrix

$P_* = [P_*; (e_{i_1})^T, (e_{i_2})^T]$, $i_1, i_2$ representing 
the index of $\bx_{\ast \ast,j}^n$, $\bx_{\ast,j}^{n+1}$ in $X^\N$. 
\State \mbox{\textbf{End While}}
\State $n_{\text{tpar}}= n_{\text{tpar}}+1$
\State \mbox{\textbf{End While}}
\State \textbf{Output:} $V_{N_{max}}$, $X^{N_{max}}_{s,j}$, $X^{N_{max}-1}_{r,j}$, $X^{N_{max}}_j$, $P_{\ast,j}$, $j=1, \cdots, n_{\text{tpar}}$. 
\end{algorithmic}
\caption{\sf AAROC: Adaptive time partitioning and adaptive enrichment strategy-based ROC algorithm}
\label{alg:aaroc}
\end{algorithm}

\section{Numerical results}
\label{sec:num}

In this section, we test the proposed AAROC method on two nonlinear time-dependent problems: the one-dimensional viscous Burgers' equation and the two-dimensional lid-driven cavity problem. A comprehensive comparison is conducted between AAROC  and AROC methods in terms of error estimator, relative error, distributions of selected parameter-time pairs, and collocation points.

To evaluate the accuracy of the AAROC and AROC methods, the average relative errors, $E_n$, are calculated for both methods using a testing set that does not overlap with the training set. The total number of testing data points is denoted by $m_{\rm{test}}$. The average relative error is defined as:
\begin{align*}
E_n = \frac{1}{m_{\rm{test}}}\sum_{i=1}^{m_{\rm{test}}} \frac{\lVert u(\cdot,\cdot;\bmu^i) -\widehat{u}_n(\cdot,\cdot;\bmu^i)\rVert_F}{\lVert u(\cdot,\cdot;\bmu^i)\rVert_F},
\end{align*}
where $\lVert \cdot \rVert_F$ denotes the Frobenius norm,
\begin{align*}
  \Vert v(\cdot,\cdot) \rVert^2_F &\coloneqq \sum_{\bx \in X^\N, t_i \in {\mathcal T}_f} v(\bx,t_i)^2.
\end{align*}

\subsection{Viscous Burgers' equation}

We first consider a one-dimensional viscous Burgers' equation \cite{Benjamin2020} with a time-dependent parametric viscosity $\nu (t)$:
\begin{align}
\begin{aligned}
&\partial_t u +u\partial_x u = \nu (t) \partial_{xx} u, ~x\in (-1,1), ~t\in (0,T]
\end{aligned}
\end{align}
where $\nu (t) =\mu \left( \sin (0.01\pi t) +2\right)$, $\mu \in (0.005,0.1)$, $T=1$, $\bmu = \mu$, and the inflow boundary condition is $u(-1,t)=2$. The parameter $\mu$ dictates the strength of the viscosity. The initial condition is
\begin{equation}
u(x,t=0) =\left\{
\begin{aligned}
2, ~~x=-1,\\
1, ~-\frac{1}{2}\leq x\leq -\frac{1}{3},\\
0, ~ else.
\end{aligned}
\right.
\end{equation}

\begin{table}[htbp]
	\centering
\setlength{\tabcolsep}{8pt}
\renewcommand{\arraystretch}{1.2}
			\begin{tabular}{|c|cccccc|}
		\hline
$p_{\rm{adap}}$ &Algorithm	&$(n_{\rm{add}},n^{\rm{incre}}_{\rm{adap}})$  &$n_{\rm{tpar}}$&$N_{\rm{adap}}$ &$E_{40}$&$E_{80}$
 \\ \hline
\multirow{2}{*}{0.1}
 &AROC&(11,5)  &1 &354&1&/ \\  \cline{2-7}
  &AAROC	&(11,5)  &2 &33&0.1021&0.0041\\  
\cline{2-7} 
\cline{2-7}
\hline
\multirow{2}{*}{0.2}
 &AROC&(11,5)   &1&53&0.1769 &0.0044\\  \cline{2-7}
 &AAROC	&(11,5)  &5&33&0.0358&0.0034\\ 
 \hline
\multirow{2}{*}{0.3}
& AROC&(11,5)  &1&244&0.0418&/ \\ \cline{2-7}
 &AAROC&(11,5)  &12&33&0.0878& 0.0036\\ \hline
	\end{tabular}
	\caption{Performance comparison for AAROC and AROC methods applied to the Viscous Burgers' equation with discontinuous initial condition. Here, $\gamma =80$, $n_0 =4$, $N_{\rm{adap}}$ is the number of total added collocation points when $N_{\max}$ basis are selected. The ``/'' denotes that the robustness indicator is not satisfied even when all candidate reduced collocation points are added.}
	\label{table:aaroc:burger}
\end{table}

\begin{figure}[thbp]
\centering
\includegraphics[scale=0.45]{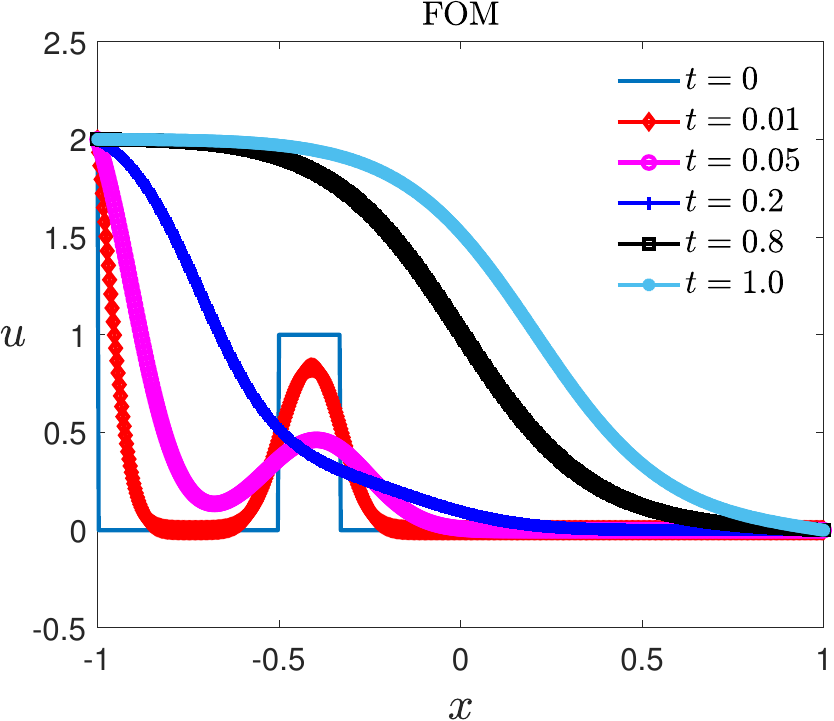}
 \includegraphics[scale=0.45]{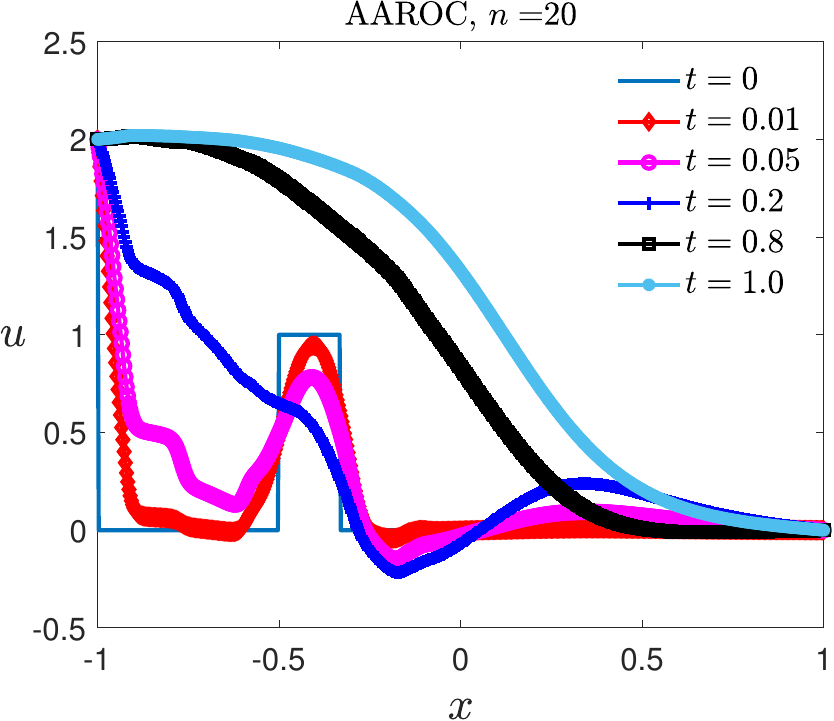}
 \includegraphics[scale=0.45]{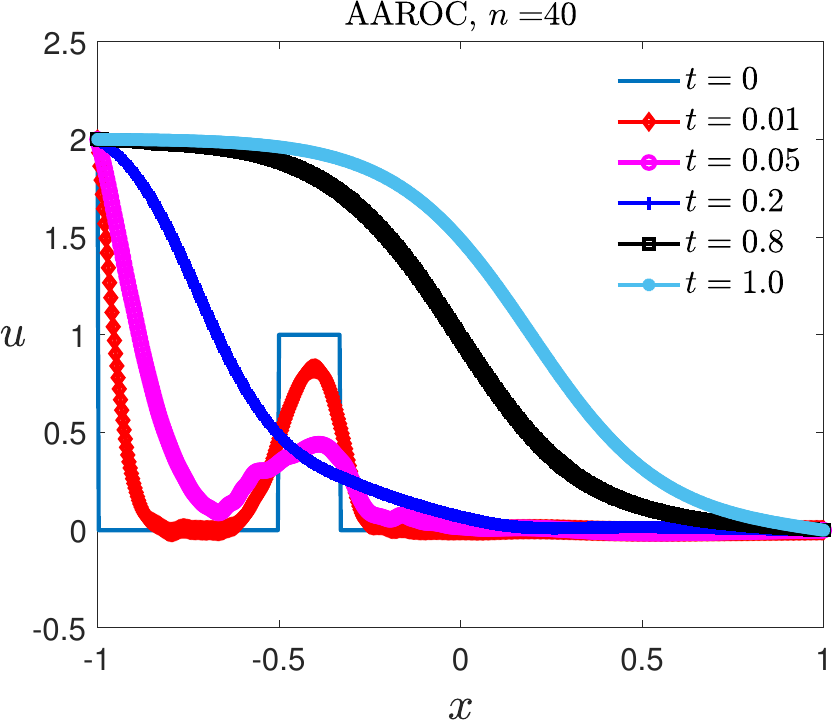}
\includegraphics[scale=0.45]{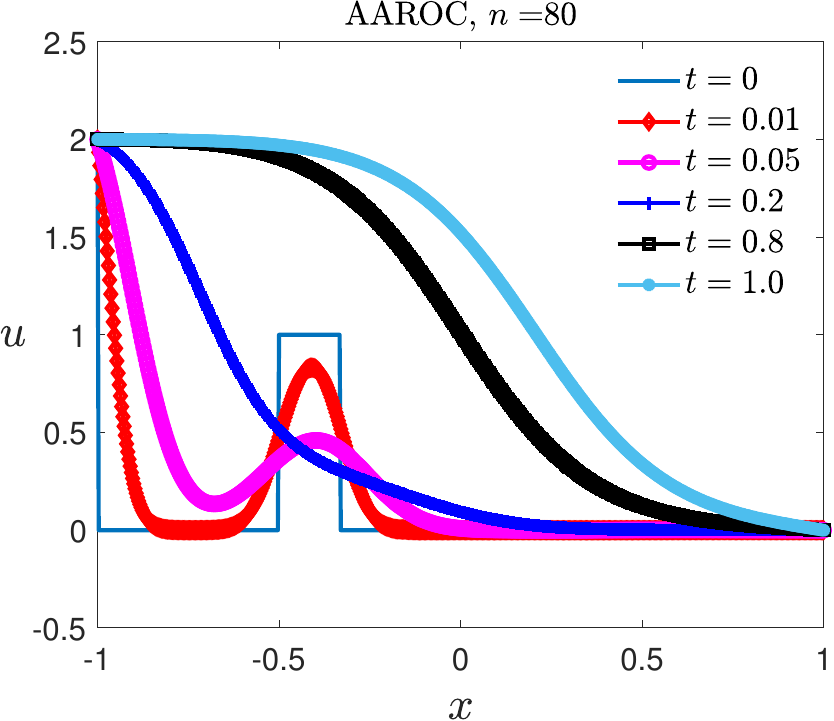}
 \caption{Results of the viscous Burgers' equation: High fidelity solution (Top Left) and AAROC approximations $\widehat{u}_n(x,t)$ at $\mu =0.0950$ with $p_{\rm{adap}} =0.2$, $n=20, 40, 80$ at $\mu =0.0950$.}
 \label{figure:burgers:fom}
\end{figure}

For the training set, 41 different values of the parameter $\mu$ are uniformly sampled on a logarithmic scale from $[0.005,0.1]$. We adopt a forward Euler time scheme with time step $\Delta t = 10^{-5}$ and the computational domain  discretized into $N_x = 600$ grids. Additionally, 7 testing parameters are uniformly sampled from the interval $[0.012, 0.095]$ , ensuring that the testing set does not overlap with the training set. Further parameters for both algorithms are detailed in Table \ref{table:aaroc:burger}. Reference high fidelity solutions for $\mu=0.095$ and the approximations of the AAROC method with various number of basis are displayed  in Figure \ref{figure:burgers:fom}. 

Our main observations are as follows. The AAROC method is more robust, accurate, and efficient than the AROC method, especially when the percentage of candidate collocation points $p_{\textrm{adap}}$ is small. 
For example, when $p_{\textrm{adap}}=0.1$, the AAROC method achieves better  accuracy than AROC with less than $10\%$ of additional reduced collocation points. It achieves a relative error of around $0.004$ when $n=80$ while AROC fail to converge robustly.

\begin{figure}[htbp]
\centering
 \includegraphics[scale=0.47]{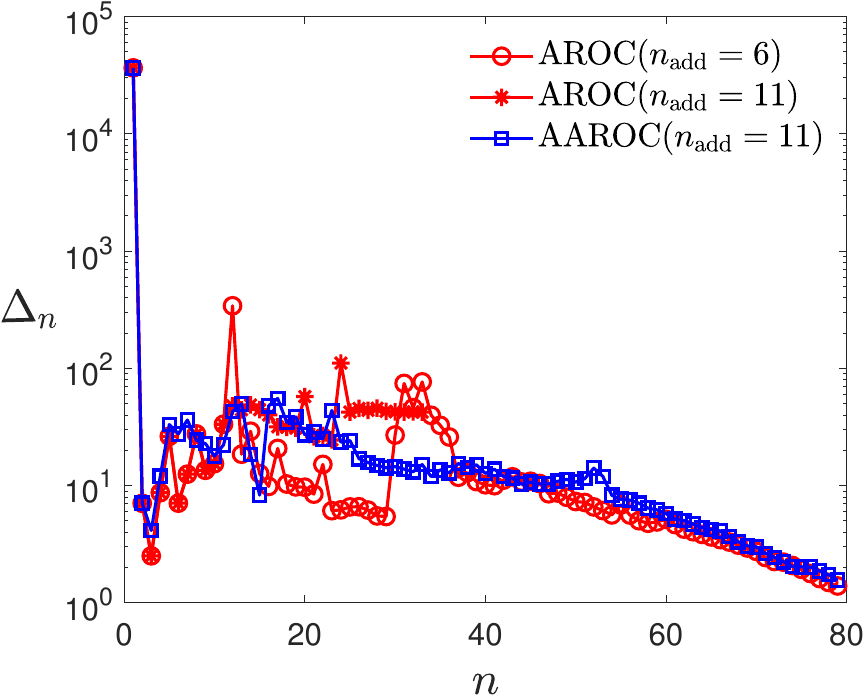}
 \includegraphics[scale=0.47]{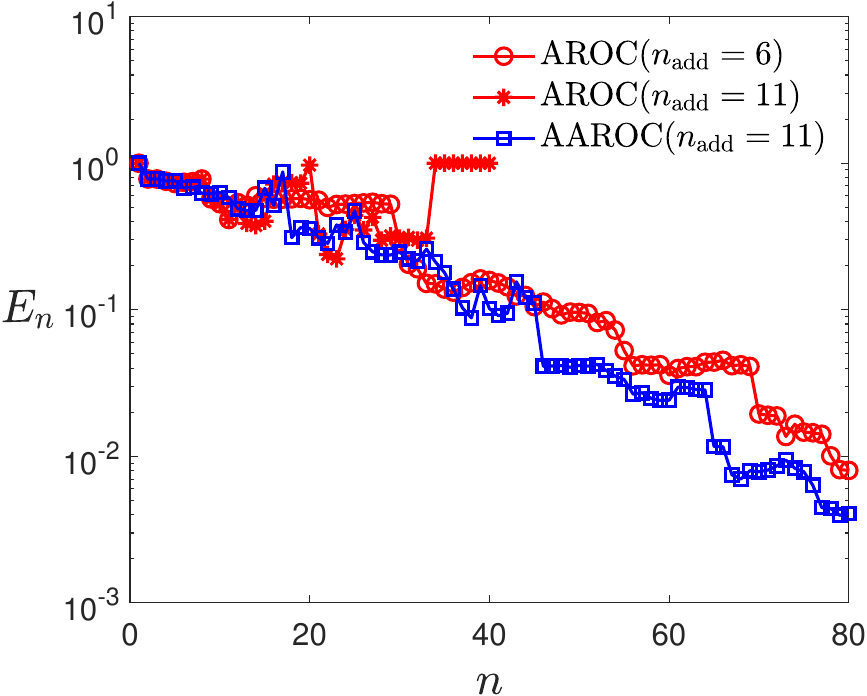}\\
\includegraphics[scale=0.47]{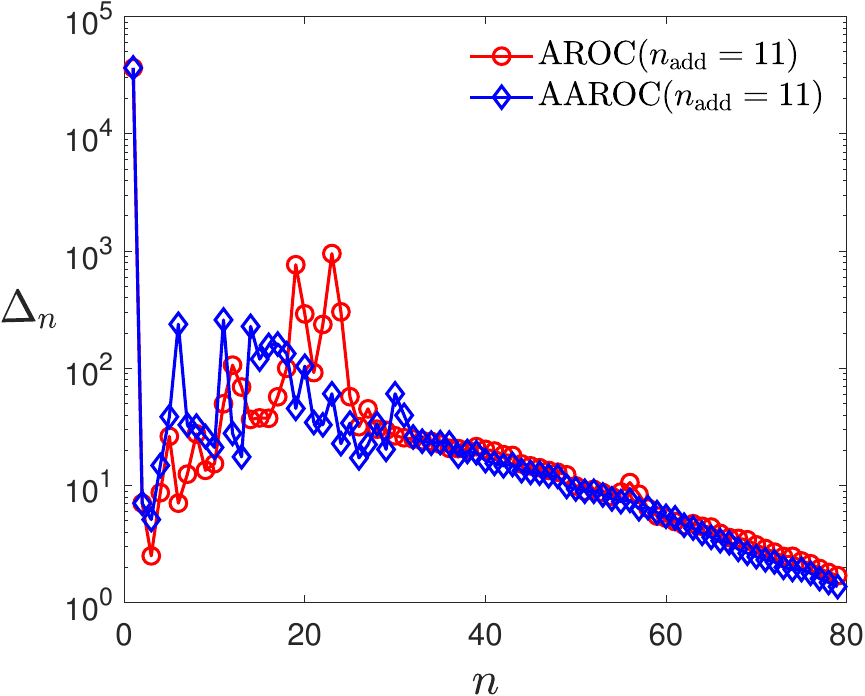}
 \includegraphics[scale=0.47]{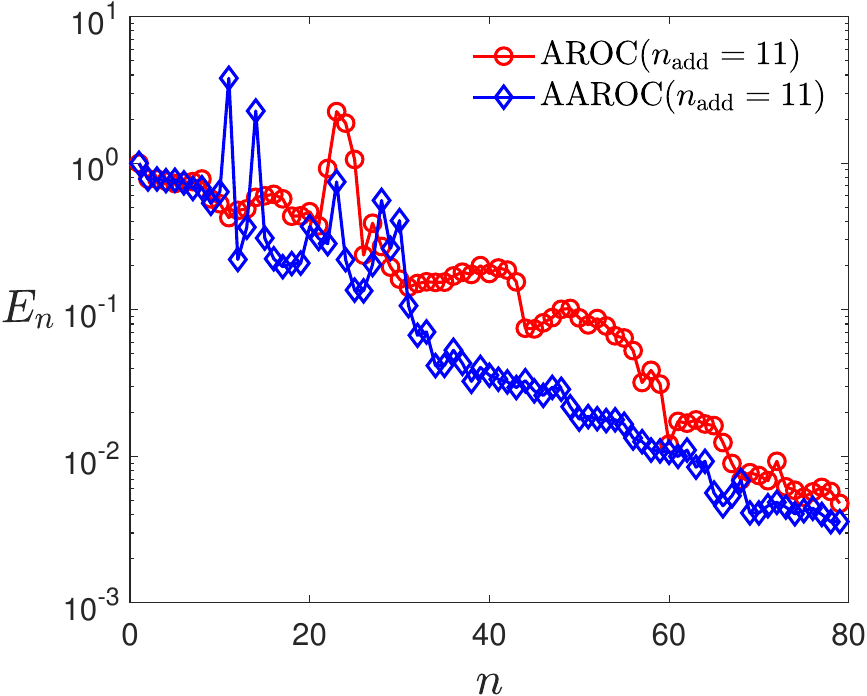}\\
\includegraphics[scale=0.47]{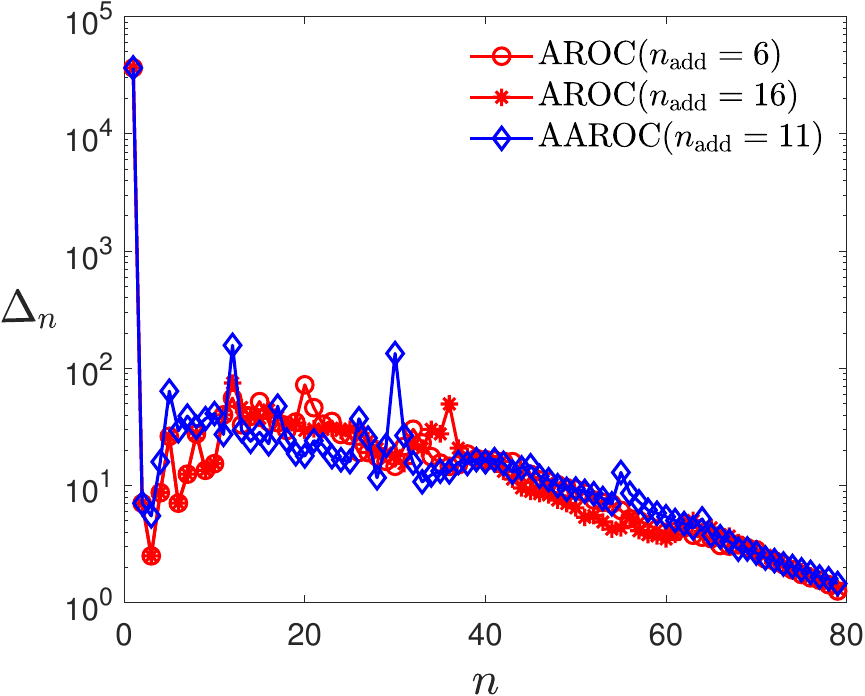}
 \includegraphics[scale=0.47]{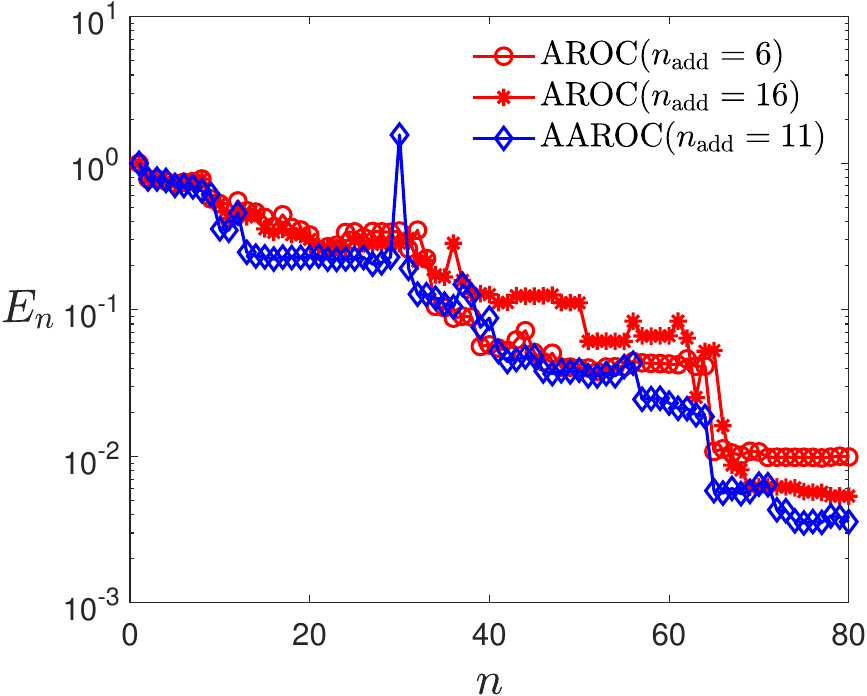}\\
	   \caption{Results of the viscous Burgers' equation: Error estimators of the offline process on the training set and relative errors of the online process on the testing set for $p_{\text{adap}}=0.1, 0.2$, and $0.3$ (Top to Bottom).}
	   \label{figure:burgers:error:all}
\end{figure}

\begin{figure}[htbp]
\centering
 \includegraphics[scale=0.47]{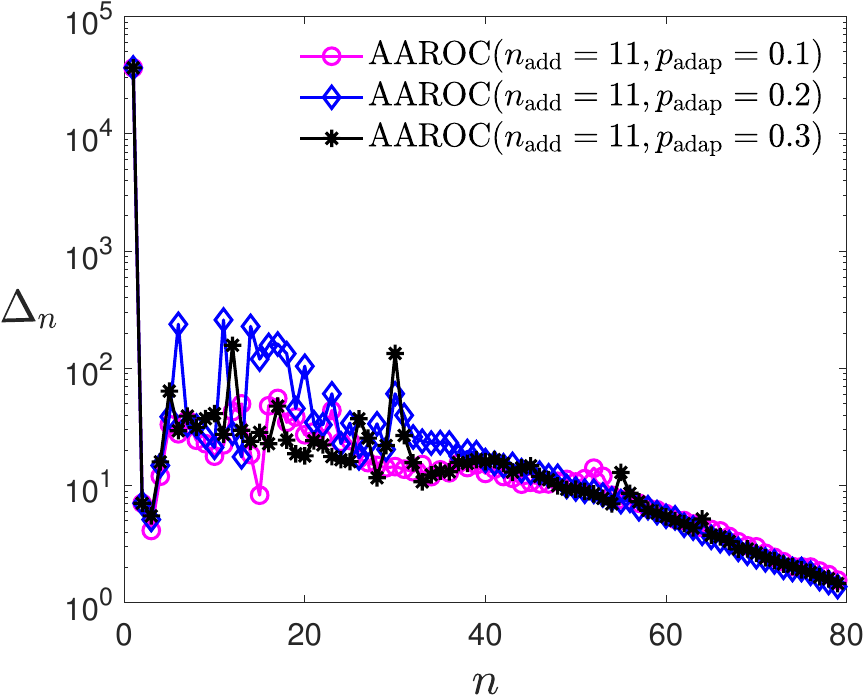}
 \includegraphics[scale=0.47]{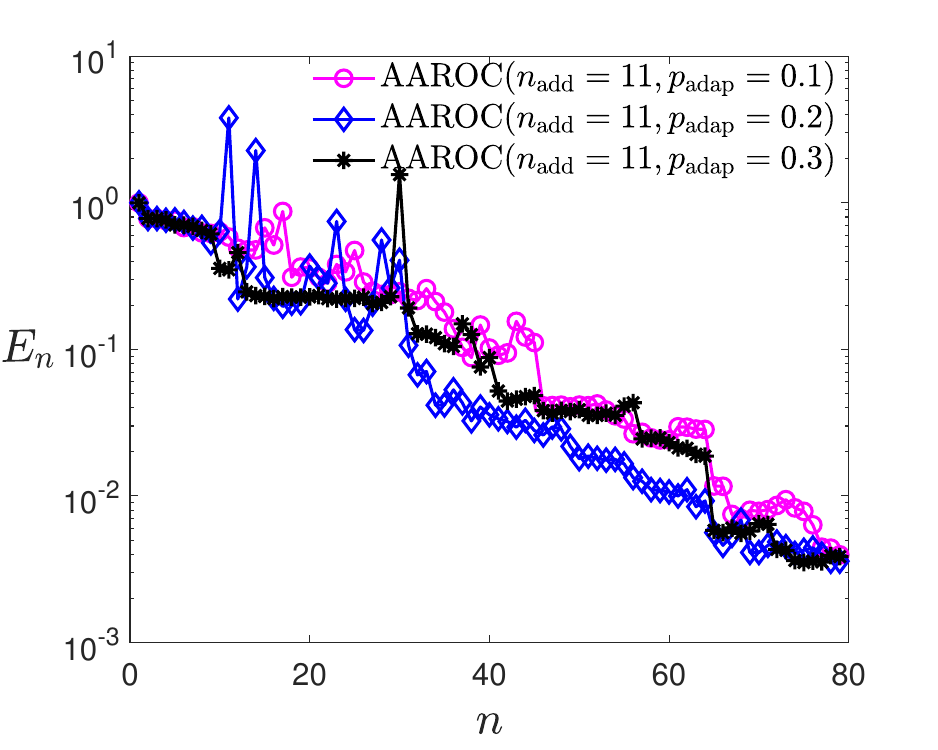}
	   \caption{Results of the viscous Burgers' equation: Error estimators of the offline process on the training set and relative errors of the online process on the testing set of the AAROC method.  }
	   \label{figure:burgers:error:n4}
\end{figure}

\begin{figure}[htbp!]
\centering
\includegraphics[scale=0.33]{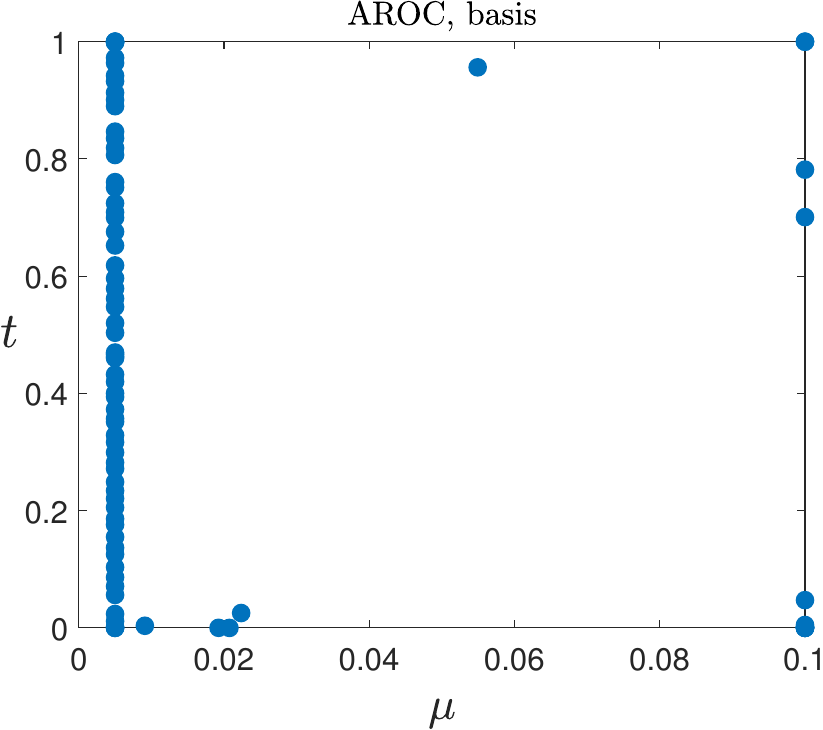}
 \includegraphics[scale=0.33]{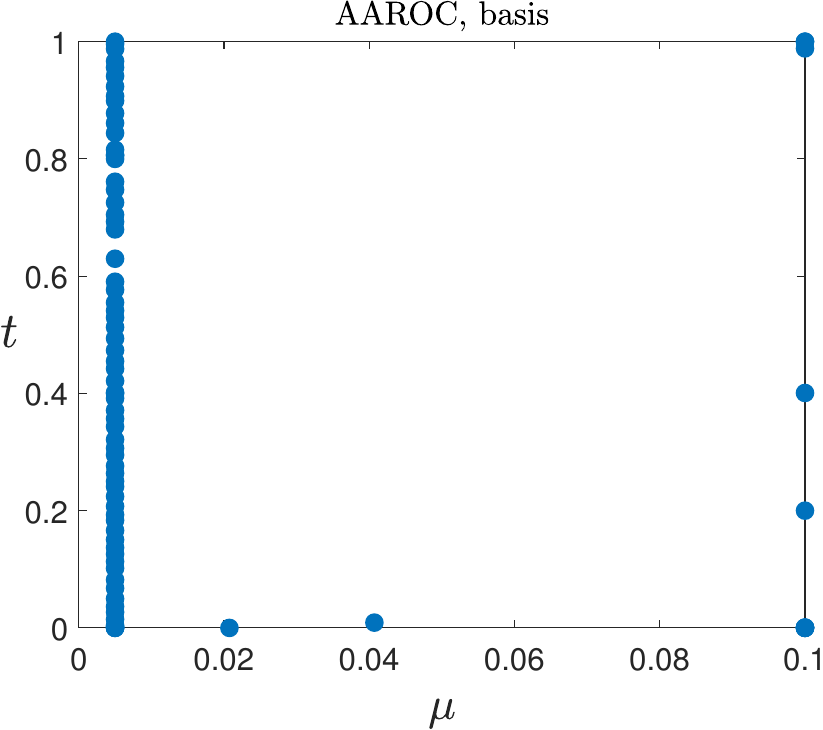}
 \includegraphics[scale=0.33]{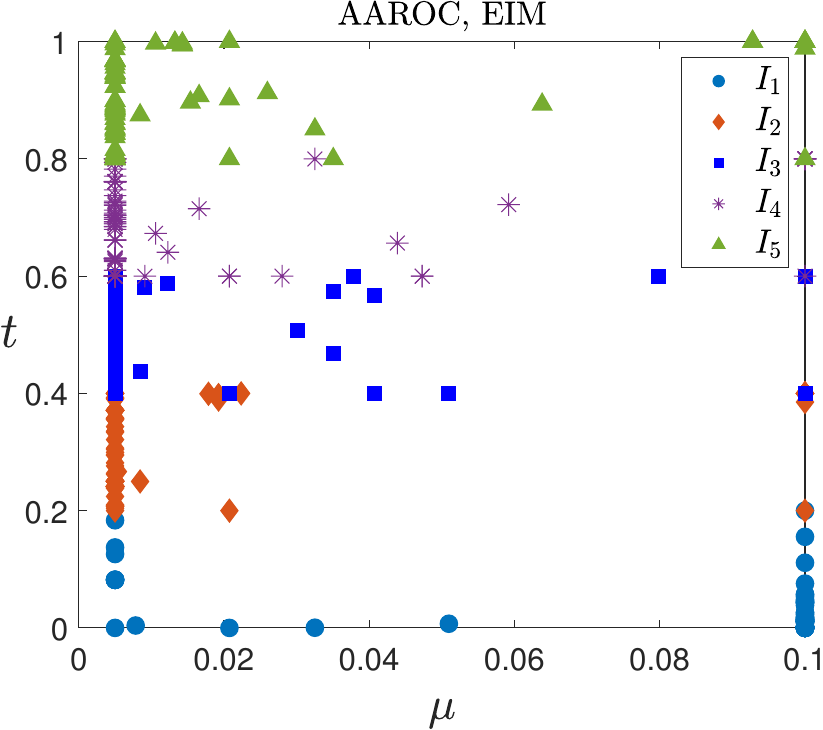}
	   \caption{Results of the viscous Burgers' equation: Distribution of selected parameter-time pairs for the AROC and AAROC basis (Left, and Middle), and that of the selected parameter-time pairs for the EIM process of the second set of collocation points of AAROC method (Right).}
	   \label{figure:burgers:collo}
\end{figure}

To provide a more intuitive understanding of the AAROC method, we consolidate numerical results from all three cases in Figures \ref{figure:burgers:error:all} and \ref{figure:burgers:error:n4}. Notably, the configuration with $p_{\text{adap}}=0.1$, and $n_{\text{add}}=31$ performs the best when the basis number exceeds 60. This outcome is expected, as it utilizes the highest number of collocation points. However, among cases with similar numbers of collocation points, the highest accuracy is achieved with $p_{\text{adap}}=0.2$, and $n_{\text{add}}=11$, as represented by the blue solid line with diamond symbols. Figure \ref{figure:burgers:collo} illustrates the effect of the time partitioning technique used in the AAROC method. The left two figures are the selected parameter-time pairs of the AROC and AAROC methods. Smaller parameter $\mu$ are more likely to be chosen in both methods. The third figure shows the distribution of selected parameter-time pair for the EIM process of the second set of collocation points. The parameter $\mu$ is sampled in a more scattered manner in later time segment.

\begin{figure}[htbp!]
\centering
 \includegraphics[scale=0.32]{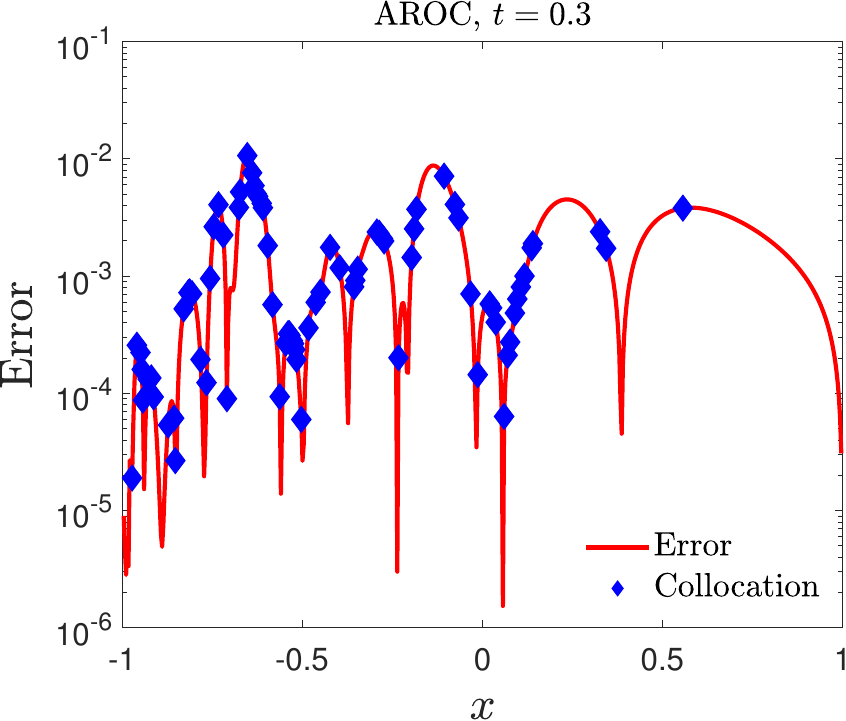}
\includegraphics[scale=0.32]{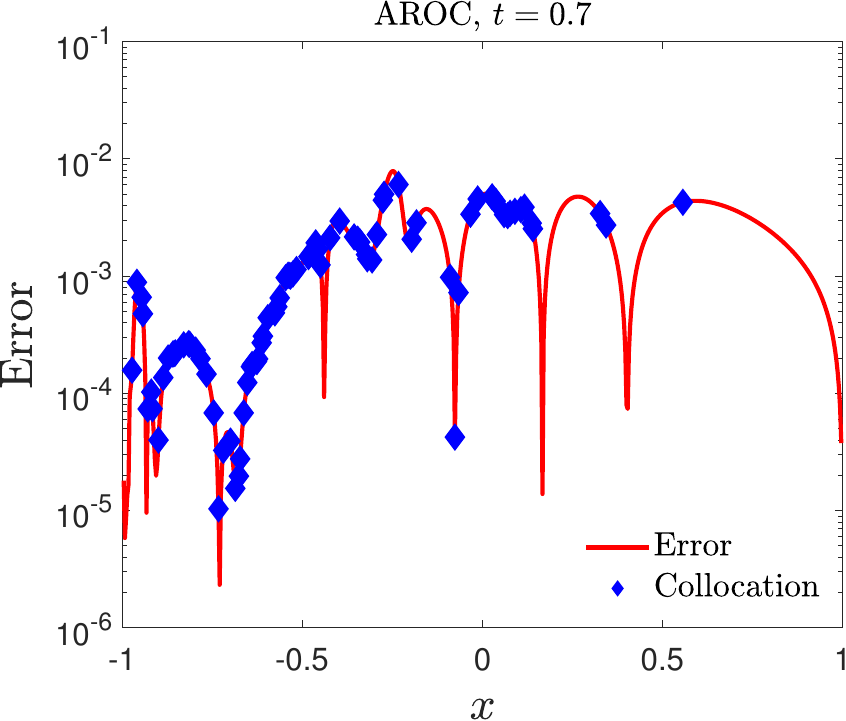}
 \includegraphics[scale=0.32]{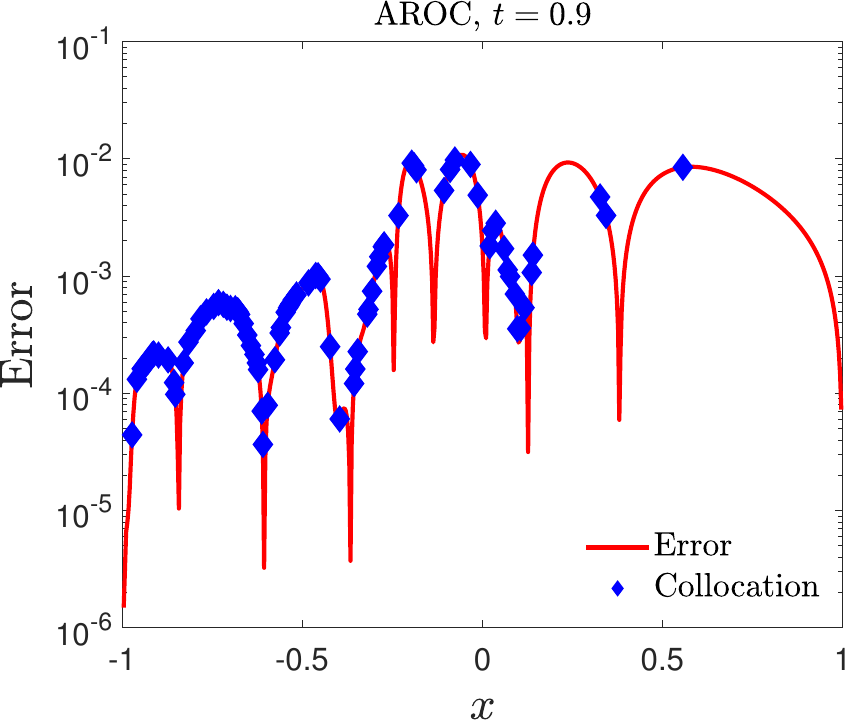}\\
  \includegraphics[scale=0.32]{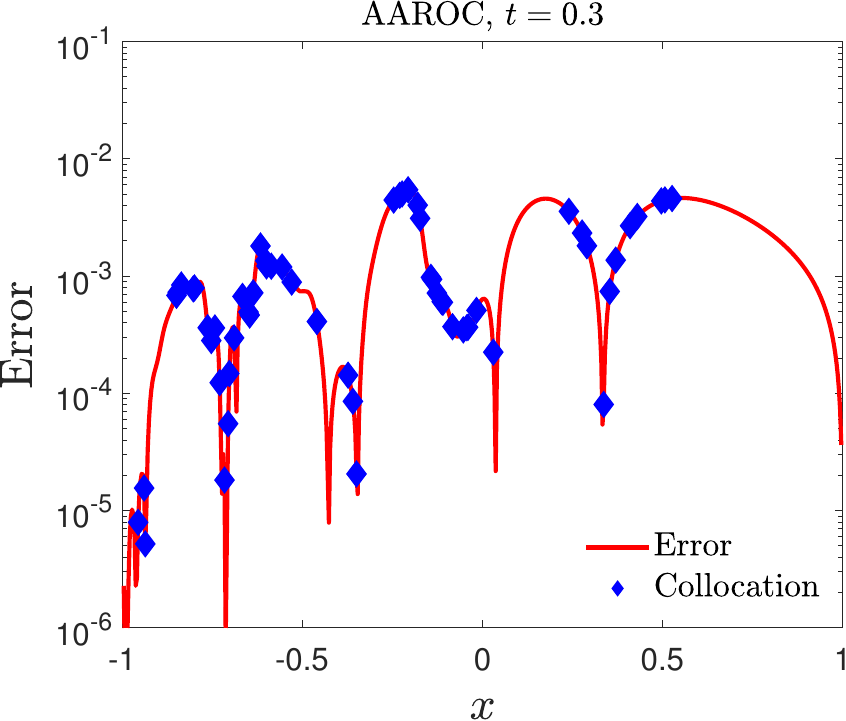}
\includegraphics[scale=0.32]{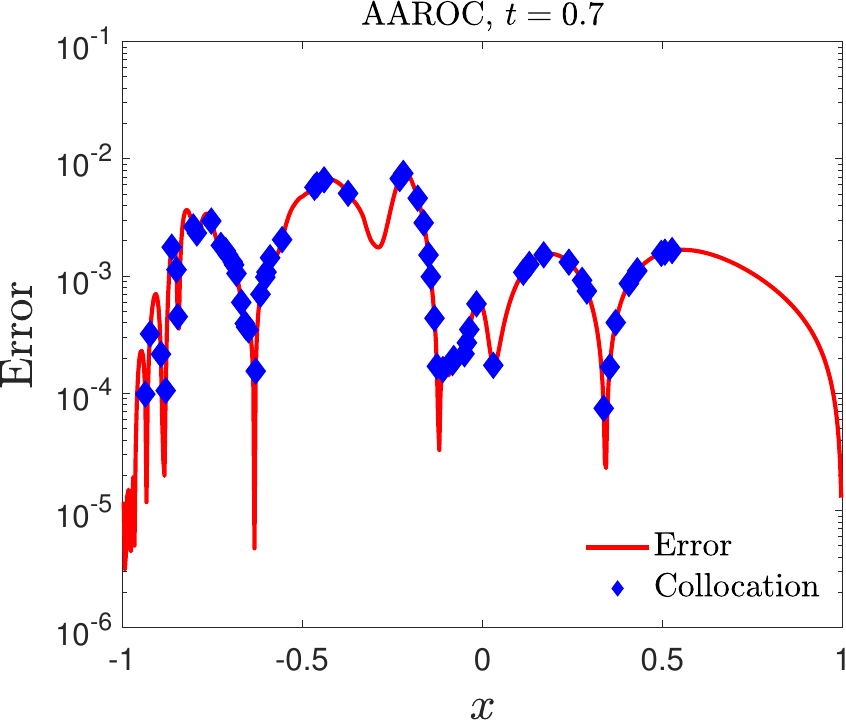}
 \includegraphics[scale=0.32]{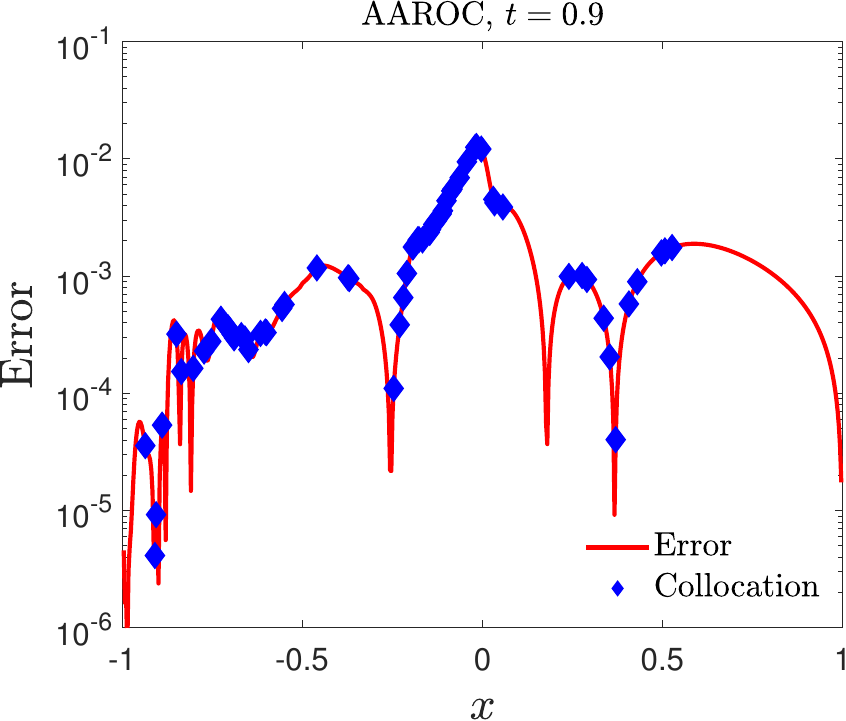}
	   \caption{Results of the viscous Burgers' equation:Relative error curves of parameter $\mu=0.012$ and corresponding unique collocation points with $n=40$ at time $t=0.3, 0.7, 0.9$ (from left to right) for the AROC (top) and AAROC (bottom) methods.}
 \label{figure:burgers:errorcollo}
\end{figure}

To illustrate the ability of the AAROC method to capture the dynamical features in the parametric problem, we present the relative error curves along with corresponding collocation points
 at time $t=0.3, 0.7, 0.9$ for parameter $\mu =0.012$ (the smallest parameter in the testing set) with $n=40$ basis in Figure \ref{figure:burgers:errorcollo}. The first row are results of AROC method, and collocation points are very sparse around $x=0.5$, where the relative error are quite large. Compared to AROC, AAROC samples more collocation points around $x=0.5$ where the relative error is larger. Moreover, the local extreme points of the error curves  are well captured by AAROC method. With the ability to adaptively add more collocation points at locations with larger error, the error of the AAROC decays faster than that of the  AROC method.

\subsection{Lid driven cavity problem}
We also consider a two-dimensional parametric lid driven cavity problem \cite{wang2012proper,stabile2018finite}. 
The equations are
 \begin{align} 
 \begin{split}            
  \frac{\partial u}{\partial t}- \frac{1}{\text{Re}}(u_{xx} + u_{yy}) +(u^2)_x +(uv)_y +\partial_x p &= 0, ~ \text{in} ~~\Omega,  \\
 \frac{\partial v}{\partial t} - \frac{1}{\text{Re}}(v_{xx} + v_{yy})+ (uv)_x + (v^2)_y +\partial_y p &=  0, ~ \text{in} ~~\Omega,  \\
  u_x + v_y &= 0, ~ \text{in} ~~ \Omega,   
  \end{split}
  \label{eq:lid}
 \end{align}
with $\Omega =[0,1] \times [0,1]$, $\text{Re} \in [10,500]$, $t\in [0,T]$, $T=35$, $\bmu =\text{Re}$. Here $u, v$, and $p$ represent the velocities in the $x$, and $y$ directions and the pressure field, respectively. Homogeneous boundary conditions are applied on all sides except for the top boundary, where the $x$-velocity $u$ is set to be $u = 1$. The spatial discretization is performed using the MAC scheme \cite{harlow1965numerical,nicolaides1992analysis,MAC_chenlong}, and backward Euler time discretization is used to obtain a high fidelity solution. The computational domain is discretized into $100 \times 100$ uniform grids. The time step is $\Delta t=0.002$. At each time step, Picard iteration is employed to solve the resulting nonlinear problem.

For the AAROC method, the training set is generated by uniformly sampling $81$ parameters from $[10, 500]$, while the testing set is generated by uniformly sampling $5$ parameters from $[10+2.2, 500-2.2]$. Other problem settings are detailed in Table \ref{table:aaroc:lid}. When $n_0=2$, AAROC samples the fewest number of collocation points while achieving the best accuracy. When $n_0=4$, the number of collocation points sampled by AAROC and AROC with $(n_{\textrm{add}},n^{\textrm{incre}}_{\rm{adap}})=(6,5)$ is comparable, but the error of the AAROC is one-fourth of that of the AROC method.

\begin{table}[htbp]
	\centering
\setlength{\tabcolsep}{8pt}
\renewcommand{\arraystretch}{1.2}
			\begin{tabular}{|c|cccccc|}
		\hline
	 $n_0$&Algorithm	&$(n_{\rm{add}},n^{\textrm{incre}}_{\rm{adap}})$  &$n_{\rm{tpar}}$&$N_{\rm{adap}}$ &$E_{20}$&$E_{45}$
 \\ \hline
 \multirow{5}{*}{2}
&AROC &(31,30)  &1 &273 &0.0840&0.0276\\  \cline{2-7}
&AROC&(11,10)  &1 &154&0.0855&0.0103\\  \cline{2-7}
&AROC &(6,5)  &1 &109&0.2080&18.6677\\  \cline{2-7}
 & AAROC  &(31,5)  &2 &93&0.0398&0.0070
 \\ 
 \hline
 \multirow{5}{*}{4}
&AROC&(31,30)  &1 &335 &0.0590&0.0152\\  \cline{2-7}
&AROC&(11,10)  &1 &163&0.0542&0.0173\\ \cline{2-7}
&AROC &(6,5)  &1 &114&0.0507&0.0155\\ \cline{2-7}
&AAROC  &(31,5)  &4 &124&0.0711&0.004
\\ \hline
	\end{tabular}
	\caption{Results of the lid driven cavity problem: Problem settings of the AAROC and AROC methods. Here, $\gamma =10$, $N_{\rm{adap}}$ is the total number of added collocation points when $N_{\max}$ basis are selected.}
	\label{table:aaroc:lid}
\end{table}

\begin{figure}[htbp!]
\centering
   \includegraphics[scale=0.42]{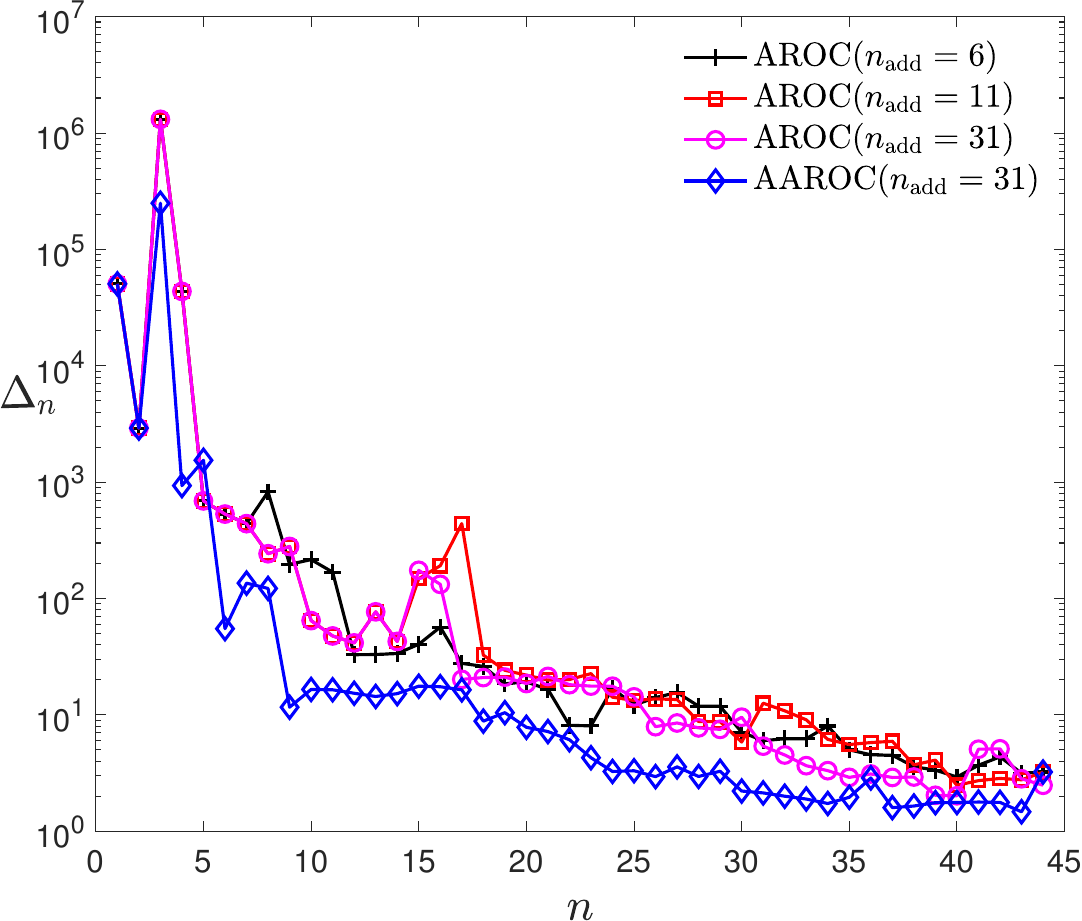}
   \includegraphics[scale=0.42]{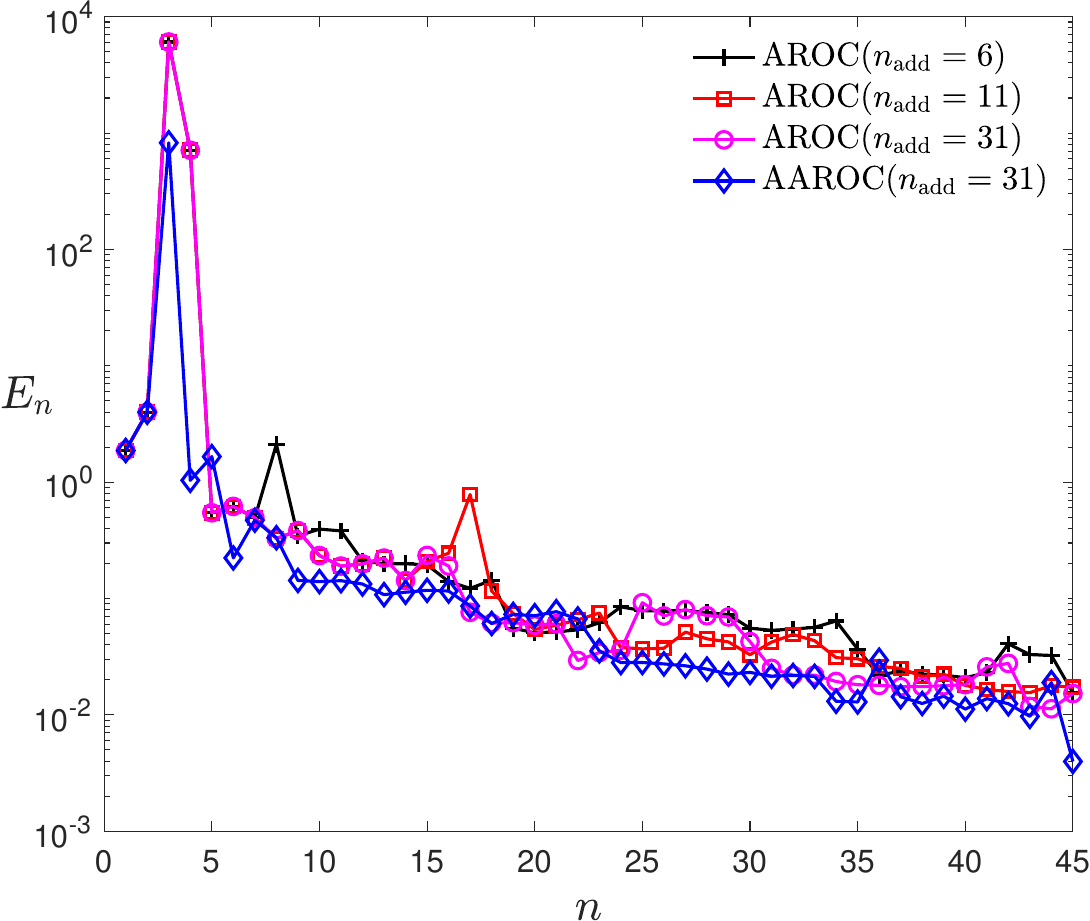}
	   \caption{Results of the lid driven cavity problem: Error estimators and relative errors of the AAROC and AROC algorithms with $n_0=4$.}
	   \label{figure:lid:error:n4}
\end{figure}
Numerical results with $n_0=4$ are presented in Figure \ref{figure:lid:error:n4}. The AAROC method adaptively adds 124 collocation points when $45$ basis are selected in the offline process. For the AROC method, three different values of $(n_{\rm{add}},n^{\textrm{incre}}_{\rm{adap}})$ are tested. The relative errors of the AAROC method are lower compared to the errors obtained with the AROC method. Noticeably, the number of additional collocation points of the AROC method with $(n_{\rm{add}}=31, n^{\textrm{incre}}_{\rm{adap}}=30)$ is nearly twice that of the AAROC method, yet the relative error in the former case is worse than that in the latter.

\begin{figure}[htbp!]
\centering
   \includegraphics[scale=0.42]{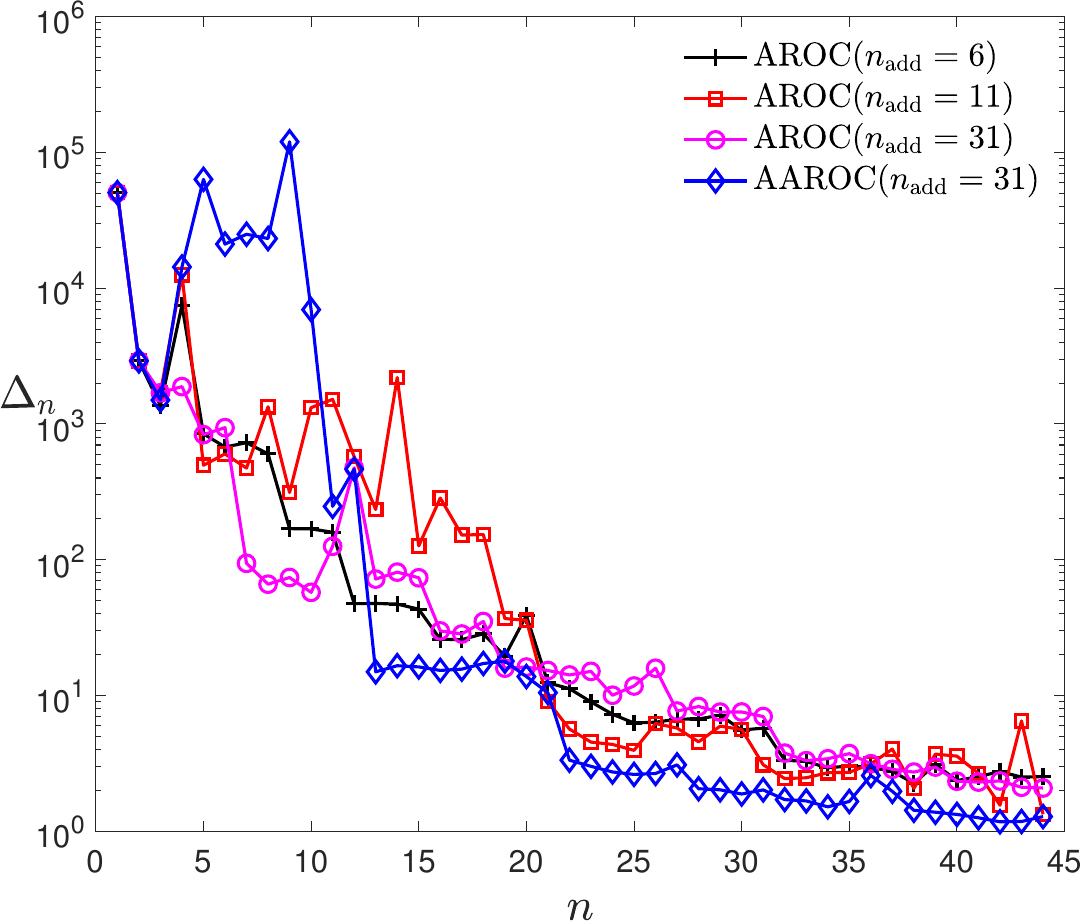}
   \includegraphics[scale=0.42]{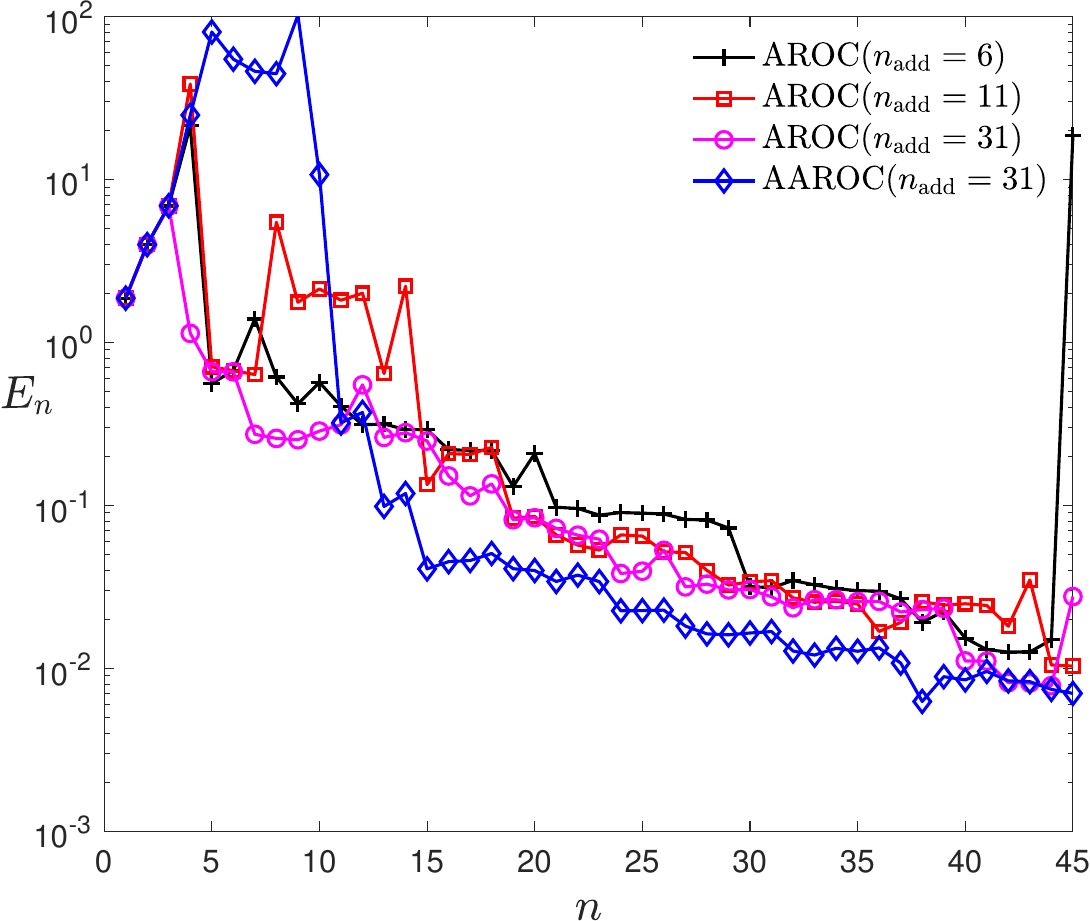}
	   \caption{Results of the lid driven cavity problem: Error estimators and relative errors of the AAROC and AROC algorithms with $n_0=2$.}
	   \label{figure:lid:error:n2}
\end{figure}

Empirical knowledge from \cite{chen2022hyper} suggests testing with $n_0=2$, which activates the adaptive enrichment strategy and adaptive time partitioning technique sooner during the offline process. The error estimators in the offline process and the relative errors of the online process are shown in Figure \ref{figure:lid:error:n2}. The AAROC method adds the fewest additional reduced collocation points while achieving the best relative error. Specifically, the blue line with diamond symbols is located at the bottom of all curves after 10 bases. Therefore, the advantages of the AAROC method in terms of both accuracy and stability are better than those of the AROC method. Additionally, the total number of added collocation points for both the AROC and AAROC methods is smaller with $n_0=2$ compared to $n_0=4$.

\begin{figure}[htbp!]
\centering
\includegraphics[scale=0.35]{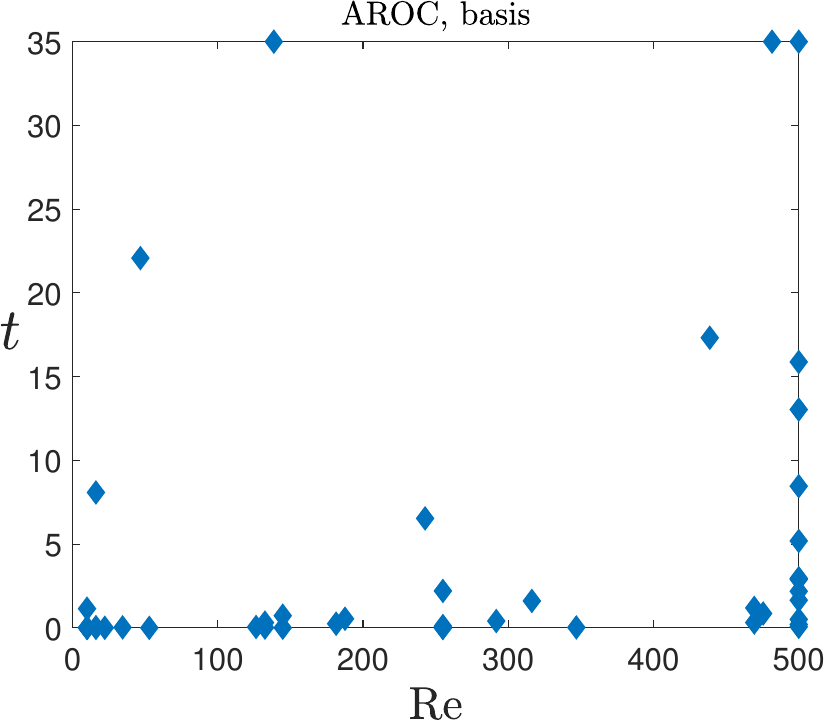}
   \includegraphics[scale=0.35]{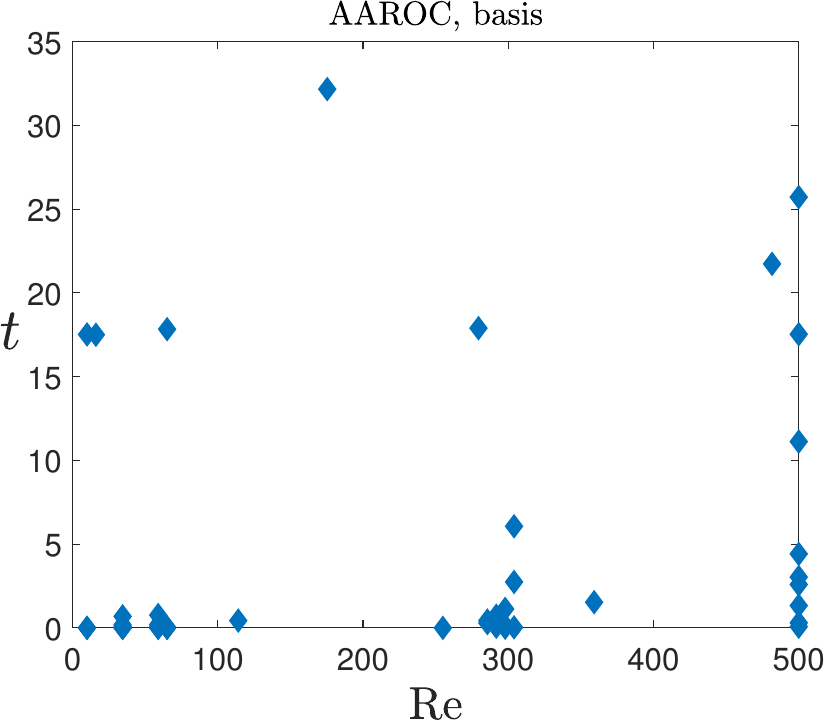}
   \includegraphics[scale=0.35]{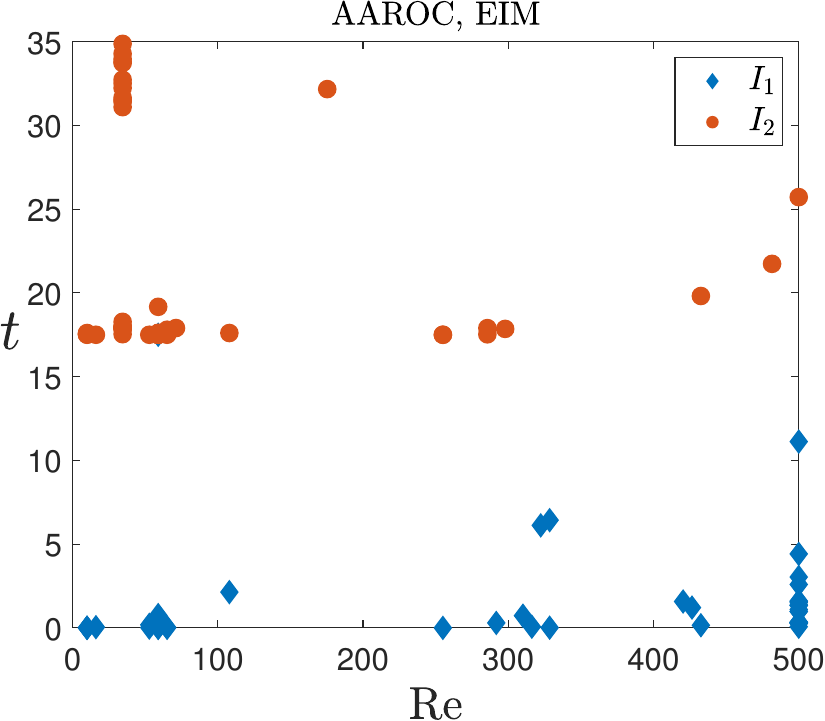}
	   \caption{Results of the lid driven cavity problem: Distributions of selected parameter-time pairs for basis of the AROC and AAROC methods, and that of the selected parameter-time pairs for the EIM process of the second set of collocation points of the AAROC method. Here $n_0=2$, $n_{\rm{add}}=11$ (AROC), $n_{\rm{add}}=31$(AAROC).
       }
	   \label{figure:lid:parapair}
\end{figure}

Figure \ref{figure:lid:parapair} shows the distributions of the selected parameter-time pairs for both basis and EIM process of the second set of collocation points.
It can be observed that the distributions of selected parameter-time pairs for the EIM process of collocation points vary across different time segments. In particular, more large time nodes have been selected in the second segment, while fewer are chosen in the first. This dynamic selection again makes contribution to the overall higher accuracy of the AAROC method compared to the AROC method.

\begin{figure}[thbp]
\centering
\includegraphics[width=0.42\textwidth]{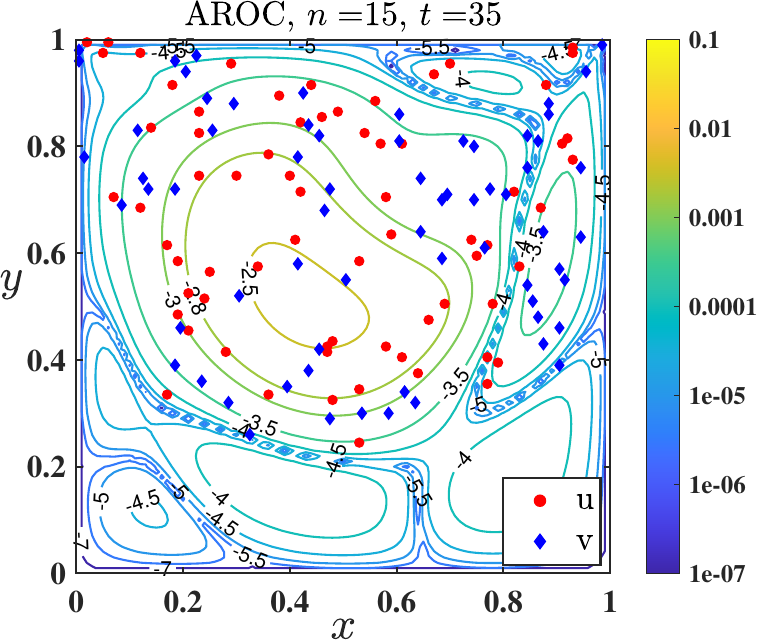}
\includegraphics[width=0.42\textwidth]{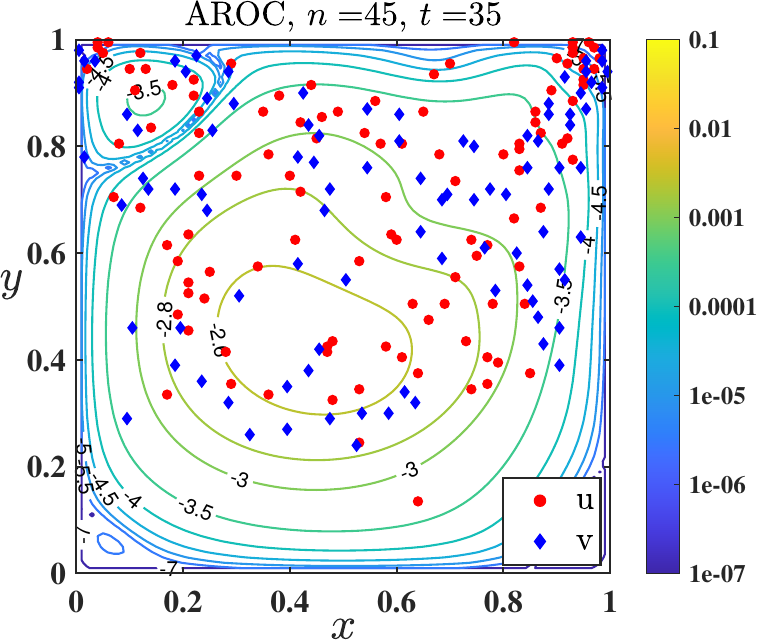}\\
\includegraphics[width=0.42\textwidth]{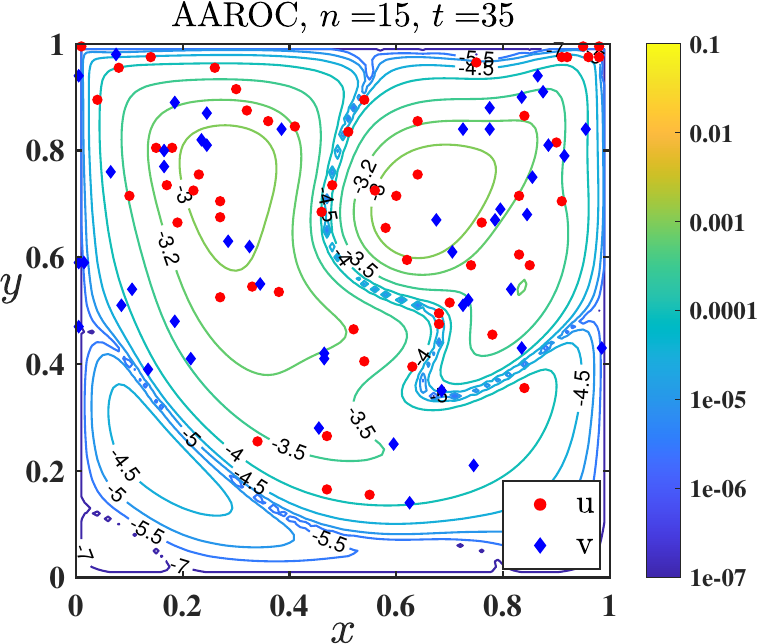}
\includegraphics[width=0.42\textwidth]{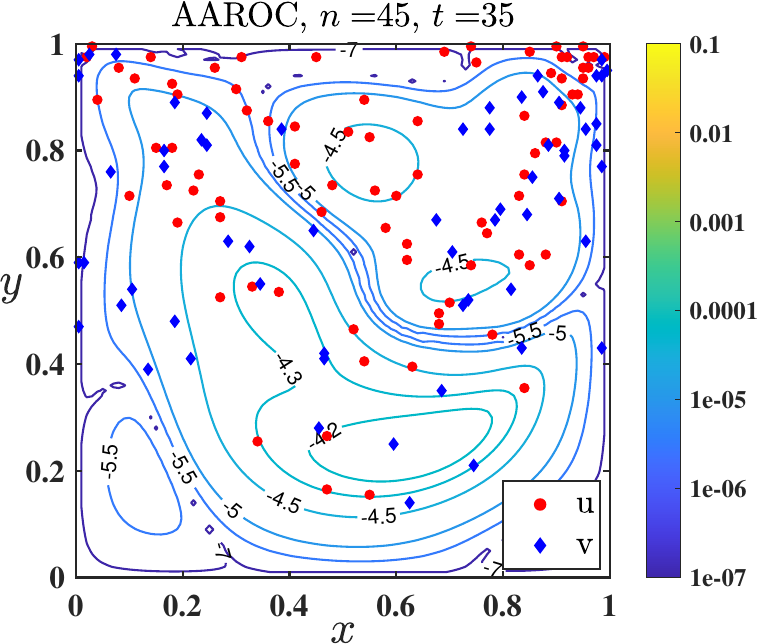}
\caption{Results of the lid driven cavity problem: Absolute error contours of streamlines of AROC (top) and AAROC (bottom) methods at final time ($t=35$), with the RB dimensions $n=15$ (left) and $n=35$ (right). Here, $\text{Re}=255$, $n_0=2$.}
\label{figure:lid:streamline} 
\end{figure}

The corresponding selected collocation points and the absolute error of streamlines of AROC and AAROC method are displayed in Figure \ref{figure:lid:streamline}. To better distinguish AROC and AAROC, we drop the common points shared by these two methods, and display only collocation points unique to each method at time $t=5$ and $t=35$. As time involves, the shape of the streamline error varies. With the help of adaptive time partitioning, AAROC method is able to more efficiently capture such variations dynamically.  With $n=15$, AAROC only samples 122 points, while AROC samples 152 points. When $n=45$, the error contours of both methods becomes smaller. The AAROC samples 182 points, while AROC samples 243 points. We note that the relative error curve of the AAROC method is always below those of the AROC method, even while sampling fewer collocation points (see Figure \ref{figure:lid:error:n2}).

\section{Conclusion}
\label{sec:conclusion}
This paper proposes a novel reduced-over collocation method with adaptive time partitioning and adaptive enrichment strategies. These two adaptive techniques are integrated into an iterative greedy framework during the offline process through the use of a robustness indicator. These approaches enhance the stability and accuracy of the AAROC method while relying on fewer collocation points, making it suitable for complex fluid dynamic problems. Further work will include extensions to a nonlinear ROM framework for better constructing surrogate models for problems with slow-decaying Kolmogorov $N$-width, and adaptive non-uniform time partitioning strategies.

\section*{Data accessibility}
The data that support the findings of this study are available from the corresponding author upon reasonable request.

\section*{Authors’ contributions}
Lijie Ji: Conceptualization, Methodology, Formal analysis, Supervision, Software, Writing - Original Draft, Writing - Review and Editing, Funding acquisition. Zhichao Peng: Conceptualization, Methodology, Formal analysis, Supervision, Writing - Original Draft, Writing - Review and Editing. Yanlai Chen: Conceptualization, Methodology, Formal analysis, Supervision, Writing - Review and Editing, Funding acquisition, Project administration. All authors have read and approved the manuscript. All authors agree to
be accountable for all aspects of the work.

\section*{Competing interests}
We declare we have no competing interests.

\bibliographystyle{abbrv}
\bibliography{main.bbl}
\end{document}